\documentclass[12pt]{article}
\usepackage[parfill]{parskip}
\usepackage[toc,page]{appendix}
\usepackage{psfrag}
\usepackage{graphicx}
\usepackage{epstopdf}

\usepackage{pinlabel}
\usepackage{enumitem}
\usepackage{amsmath}
\usepackage{amssymb}
\usepackage{amscd}
\usepackage{picinpar}
\usepackage{times}
\usepackage{pb-diagram}
\usepackage{graphicx}
\usepackage{wrapfig}
\usepackage{xspace}
\usepackage{hyperref}
\usepackage{xcolor}
\usepackage{url}
\usepackage{float}
\usepackage{algorithm}
\usepackage[noend]{algpseudocode}

\newtheorem{theorem}{Theorem}[section]

\newtheorem{proposition}[theorem]{Proposition}

\newtheorem{definition}[theorem]{Definition}

\newtheorem{problem}[theorem]{Problem}

\begin{document}

\title{Computing Harmonic Maps and Conformal Maps on Point Clouds}
\date{}

\author{Tianqi Wu
\thanks{Center of Mathematical Sciences and Applications, Harvard University, Cambridge, MA 02138. Email:
tianqi@cmsa.fas.harvard.edu
}
\and
Shing-Tung Yau
\thanks{Department of Mathematics, Harvard University, Cambridge, MA 02138. Email:
yau@math.harvard.edu
}
}



\maketitle

\begin{abstract}
We propose a novel meshless method to compute harmonic maps and conformal maps for surfaces embedded in the Euclidean 3-space, using point cloud data only. Given a surface, or a point cloud approximation, we simply use the standard cubic lattice to approximate its $\epsilon$-neighborhood. Then the harmonic map of the surface can be approximated by discrete harmonic maps on lattices. The conformal map, or the surface uniformization, is achieved by minimizing the Dirichlet energy of the harmonic map while deforming the target surface of constant curvature. We propose algorithms and numerical examples for closed surfaces and topological disks.
\end{abstract}
\section{Introduction}

Roughly speaking, a map between two surfaces is called \emph{conformal} if it preserves angles, and is called \emph{harmonic} if it minimizes the stretching energy. Computing harmonic maps and conformal maps has a wide range of applications, such as surface matching, surface parameterization, shape analysis and so on. See \cite{koehl2013automatic}\cite{gu2004genus}\cite{boyer2011algorithms}\cite{hong2006conformal}\cite{lui2010optimized}\cite{lipman2009mobius}\cite{kim2011blended}\cite{soliman2018optimal}\cite{maron2017convolutional}\cite{ben2018multi}\cite{kittler2018conformal} for examples of applications of conformal maps, and \cite{li2009surface}\cite{pinkall1993computing}\cite{joshi2007surface}\cite{li2001moving}\cite{wang2008high} for examples of applications of harmonic maps.

Existing methods for computing harmonic maps and conformal maps mostly rely on the triangle mesh approximation of a surface. However, it is often much easier to get point cloud data, rather than the triangle mesh data. Contemporary 3D scanners can easily provide 3D point cloud data sampled from the surfaces of solid objects, but sometimes it is inconvenient to generate meshes upon point clouds. Since point clouds data do not contain information about the connectivity, a lot of algorithms, which were well-established on meshes, cannot be extended to point clouds.

This paper introduces a novel meshless method for computing harmonic maps and conformal maps
for a surface in the Euclidean 3-dimensional space. The basic idea is to use a dense 3-dimensional lattice to approximate the $\epsilon$-neighborhood of the surface, and then compute the discrete harmonic map from the lattice to the target surface, by minimizing the Dirichlet energy (i.e., the stretching energy). Conformal maps are computed by minimizing the Dirichlet energy of the harmonic maps, as we deform the target surface of constant Gaussian curvature.
In this paper, we focus on harmonic diffeomorphisms and conformal diffeomorphisms to surfaces of constant curvature $\pm1$ or $0$. These maps are particularly useful for global surface parameterizations. More specifically, we propose algorithms and numerical examples for (1) maps to the unit sphere, and (2) maps to flat rectangles, and (3) maps to flat tori, and (4) maps to closed hyperbolic surfaces.

\subsection{Previous Works}
Comparing to methods for triangle meshes, there are much fewer works on meshless methods of computing conformal and harmonic maps, especially for higher genus surfaces.
Guo et al. \cite{guo2006meshless} computed global conformal parameterizations of surfaces by computing holomorphic 1-forms on point clouds.
Li et al. \cite{li2009meshless} computed harmonic volumetric maps by grids discretizations.
Meng-Lui \cite{meng2015theory} developed the theory of computational quasiconformal geometry on point clouds.
Using approximations of the differential operators on point clouds,
Liang et al. \cite{liang2013solving} and Choi et al. \cite{choi2016spherical} computed the spherical conformal parameterizations of genus-0 closed surfaces,
and
Meng et al. \cite{meng2016tempo} computed quasiconformal maps on topological disks.
Li-Shi-Sun \cite{li2017compute} computed quasiconformal maps from surfaces to planar domain, using the so-called point integral method for discretizing integral equations for point clouds.

There is an extensive literature on computing conformal maps for triangle meshes.
Gu-Yau \cite{gu2002computing}\cite{gu2003global} developed the method of computing conformal structures of surfaces by computing the discrete holomorphic one-forms.
Pinkall-Polthier \cite{pinkall1993computing}  proposed a method of conformal parameterization by computing a pair of conjugate harmonic functions.  L\'evy et al. \cite{levy2002least} and  Lipman \cite{lipman2012bounded} and Lui et al. \cite{lui2014teichmuller} computed conformal or quasiconformal maps by minimizing or controlling the conformal distortion.
There is also a big family of methods based on various notions of discrete conformality for triangle meshes, such as circle patterns \cite{rodin1987convergence}\cite{chow2003combinatorial}\cite{hurdal2009discrete}\cite{kharevych2006discrete}, and inversive distances \cite{bowers2004uniformizing}\cite{bowers2003planar}, and vertex scalings  \cite{luo2004combinatorial}\cite{bobenko2015discrete}\cite{kharevych2006discrete}, and modified vertex scalings \cite{gu2018discrete} \cite{gu2018discrete2} \cite{sun2015discrete}\cite{springborn2017hyperbolic} allowing diagonal switches. Some related convergence results for discrete conformality can be found in \cite{rodin1987convergence}\cite{he1998thec}\cite{gu2019convergence}\cite{Luo2016discrete}\cite{bucking2018c}, and other mathematical analysis can be found in \cite{de1991principe}\cite{guo2011local}\cite{glickenstein2005combinatorial}\cite{glickenstein2005maximum}\cite{glickenstein2011discrete}.
Other works on computing conformal maps on triangle meshes include
\cite{bishop2010conformal}\cite{driscoll1998numerical}\cite{stephenson2005introduction}\cite{trefethen1998schwarz}\cite{levy2002least}\cite{desbrun2002intrinsic}\cite{yueh2017efficient}\cite{dym2019linear}\cite{crane2013robust}\cite{wu2015rigidity}.

For computing harmonic maps on triangle meshes, Gaster-Loustau-Monsaingeon
\cite{gaster2018computing}\cite{gaster2019computing} give detailed discussions on the mathematical analysis and algorithms, and produced a computer software.
Notions of discrete Dirichlet energy have been discussed and used to compute discrete harmonic maps, not only in  Gaster-Loustau-Monsaingeon's work \cite{gaster2018computing}\cite{gaster2019computing}, but also extensively in other works such as \cite{de1991comment}\cite{pinkall1993computing}\cite{korevaar1993sobolev}\cite{jost2007harmonic}\cite{wang2000generalized}\cite{izeki2005combinatorial}\cite{eells2001harmonic}\cite{fuglede2008homotopy}\cite{hass2015simplicial}.


\subsection{Contribution}
To the best of the authors' knowledge, our method is the first meshless method in computing harmonic and conformal maps to surfaces of genus greater than 1.

Existing meshless methods always use approximating differential operators on point clouds, to compute harmonic or conformal maps. Our method
provides a novel type of approach, by jumping out of the point clouds to cubic lattices. It also has the following nice properties.
\begin{itemize}
\item The idea of lattice approximation is conceptually simple, and could be implemented for any topological types of surfaces.

\item One can simply tune the density of the lattices to get satisfactory accuracy, within the ability of the computing powers.

\item Numerical experiments indicate that sparse lattices already work well.

\item The lattice structure should be suitable for parallel computation, and multi-grid methods.
\end{itemize}

We proposed algorithms of computing harmonic maps and conformal maps for any closed surfaces and topological disks, embedded in $\mathbb R^3$, to constant curvature surfaces.

\subsection{Notations}
Given a point, or equivalently a vector $x$ in $\mathbb R^3$, denote
$$
|x|_2=\sqrt{x_1^2+x_2^2+x_3^2}
$$
as the $l^2$-norm of $x$.

Given two points $x,y\in\mathbb R^3$, denote
$$
d(x,y)=|x-y|_2
$$
as the distance between $x$ and $y$.

Given a subset $A\subset\mathbb R^3$ and a point $x\in\mathbb R^3$, denote
$$
d(x,A)=\inf_{y\in A}d(x,y)
$$
as the distance from $x$ to $A$.

Given a subset $A\subset\mathbb R^3$ and $r>0$, denote
$$
B(A,r)=\{x\in\mathbb R^3:\exists y\in A, s.t. ~|x-y|_2<r\}
$$
as the $r$-neighborhood of the subset $A$.

Given a closed surface $M$ embedded in $\mathbb R^3$, we say a finite subset $P\subset\mathbb R^3$ is a $\delta$-\emph{point cloud} of $M$ if
$$
P\subset B(M,\delta) \quad\text{and}\quad M\subset B(P,\delta).
$$

If $G=(V,E)$ is an undirected connected simple graph, and $w\in\mathbb R^{E}_{>0}$ is an edge weight, then $L=L_{G,w}$ denotes the discrete Laplacian such that
for any $f:V\rightarrow\mathbb R^k$,
$$
(Lf)(i)=\sum_{j:ij\in E}w_{ij}(f(j)-f(i)).
$$
\subsection{Organization of the Paper and Acknowledgement}

The remaining of the paper is organized as follows. We review the basic mathematical theory in Section 2, and then introduce the lattice approximation of a surface in Section 3.
The algorithms, as well as numerical examples, are given in Section 4, 5, 6, 7, for the cases of
 spheres, rectangles, flat tori, and hyperbolic surfaces respectively.

This work is partially supported by Center of Mathematical Sciences and Applications at Harvard University, and Yau Mathematical Sciences Center at Tsinghua University, and NSF 1760471.

\section{Mathematical Background}
\subsection{Conformal Maps and the Uniformizations}
We give a formal definition of a conformal map as follows.
\begin{definition}[Conformal Map]
A diffeomorphsim $f$ between two Riemannian surfaces $(M,g)$ and $(N,h)$ is called a \emph{conformal map} if $f^*h=\lambda^2g$ for some smooth positive function $\lambda$, or equivalently, the local angles are preserved under $f$.
\end{definition}
In this paper we always assume that a conformal map is a diffeomorphism between two surfaces.
By the celebrated uniformization theorem, it is well known that any closed orientable surface is conformally equivalent to a surface of constant Gaussian curvature.

\begin{theorem}[Uniformization for Closed Surfaces]
Given a closed orientable Riemannian surface $(M,g)$, there exists a Riemannian surface $(N,h)$, with constant Gaussian curvature $K$, and a conformal map $f:M\rightarrow N$ such that
\begin{enumerate}[label=(\alph*)]
\item $K=1$ if $M$ has genus $0$ and

\item $K=0$ if $M$ has genus $1$, and

\item $K=-1$ if $M$ has genus $>1$.
\end{enumerate}
\end{theorem}

The Riemannian surface $(N,h)$ above is called a \emph{uniformization} of $(M,g)$ and could be viewed as a canonical representative in the conformal equivalence class.
A closed surface $(N,h)$ of constant curvature $-1$ or $0$ can be naturally represented as a planar domain, after a few cuts.
\begin{enumerate}[label=(\alph*)]
\item If $(N,h)$ has constant curvature 0, then it is isometric to the flat complex plane $\mathbb C$ modulo a lattice $\mathbb Z\oplus\mathbb Z\tau$ where $\tau\in\mathbb C$ and $Im(\tau)>0$. If we properly cut two loops on $N$, then $(N,h)$ can be displayed isometrically as a planar domain in $\mathbb C$.

\item if $(N,h)$ has constant curvature -1, then it is isometric to the hyperbolic plane $\mathbb H^2$ modulo a discrete subgroup $\Gamma$ of the orientation-preserving isometry group $Isom^+(\mathbb H^2)$. In the Poincar\'e disk model, the hyperbolic plane $\mathbb H^2$ is identified as the unit disk $\{z\in\mathbb C:|z|<1\}$ with the complete Riemannian metric
$$
\frac{4|dz|^2}{(1-|z|^2)^2}=\frac{4(dx^2+dy^2)}{(1-x^2-y^2)^2}.
$$
Similar to the flat case, if we properly cut a few loops, $(N,h)$ can be displayed, isometrically in the hyperbolic sense, as a domain in the unit disk $\{|z|<1\}$.
\end{enumerate}
So for genus 0 surfaces, there is essentially a unique conformal equivalence class. For genus 1 surfaces, the conformal equivalence classes can be parameterized by a complex number $\tau$ with $Im(\tau)>0$. For a surface $M$ with genus $>1$, the conformal equivalence classes can be parameterized by the generators of the discrete subgroup $\Gamma$ of $Isom^+(\mathbb H^2)$.
As shown above, a surface of constant curvature $\pm1$ or $0$ can always been identified as a planar domain, by a stereographic projection or cutting loops. So computing diffeomorphisms to such surfaces gives rise to global parameterizations and surface flattening, which have fundamental applications in computational graphics.

For nonclosed surfaces, the most important case is the topological disk. Assume $(M,g)$ is a Riemannian surface homeomorphic to a closed disk, then again by the uniformization theorem it is conformally equivalent to a closed unit disk. However, for our convenience in utilizing the method of harmonic maps, we consider conformal maps and harmonic maps to rectangles instead of disks. The following uniformization-type theorem for rectangles is well-known.

\begin{theorem}
Assume $(M,g)$ is a Riemannian surface homeomorphic to a closed disk. Given 4 ordered points $A_1,A_2,A_3,A_4$ on $\partial M$, there exists unique $a\in\mathbb R_{>0}$ and a conformal map
$$
f:M\rightarrow [0,a^{-1}]\times[0,a]
$$
such that
$$
f(A_1)=(0,0),\quad
f(A_2)=(a^{-1},0),\quad
f(A_3)=(a^{-1},a),\quad
f(A_4)=(0,a).
$$
\end{theorem}

\subsection{Harmonic Maps}
Let $(M,g)$ and $(N,h)$ be two smooth Riemannian surfaces. The \emph{Dirichlet energy}, i.e., the stretching energy of a smooth map $f:M\rightarrow N$ is formally defined as
$$
\mathcal E(f)=\frac{1}{2}\int_M\|df\|^2dv_g
$$
where $dv_g$ denotes the volume element on $(M,g)$, and $\|df\|$ is the norm of the differential of $f$, with respect to the induced metric on $T^*M\otimes f^*(TN)$.

\begin{definition}
A smooth map $f:M\rightarrow N$ is called \emph{harmonic} if it is a critical point of the Dirichlet Energy.
\end{definition}
Smooth harmonic maps have been extensively studied. See \cite{schoen1979existence}\cite{eells1995two} for examples.
Here we focus on harmonic diffeomorphisms to surfaces of constant curvature $\pm1$ or $0$, and harmonic diffeomorphisms to rectangles with given boundary correspondences.
Some relevant well-known facts are summarized as the following 2 theorems.
\begin{theorem}
Assume $M,N$ are two closed orientable Riemannian surfaces, and $N$ has constant curvature $K=\pm1$ or 0, and $f$ is a diffeomorphism from $M$ to $N$.
\begin{enumerate}[label=(\alph*)]

\item $
\mathcal E(f)\geq Area(N),
$
and the equality holds if and only if $f$ is conformal.

\item If $f$ is harmonic, $f$ minimizes the Dirichlet energy in its homotopy class.

\item If $f$ is conformal, then $f$ is harmonic.

\item If $f$ is harmonic and $N$ is a unit sphere, then $f$ is conformal.

    \item $f$ is always homotopic to a harmonic diffeomorphism, which is (i) unique up to a M\"obius transformation if $K=1$, and (ii) unique up to a translation if $K=0$, and (iii) unique if $K=-1$.
\end{enumerate}
\end{theorem}

\begin{theorem}
Assume
\begin{enumerate}[label=(\roman*)]
\item $(M,g)$ is a Riemannian surface homeomorphic to a closed disk, and

\item $A_1,A_2,A_3,A_4$ are 4 ordered points on $\partial M$, and

\item $a\in\mathbb R_{>0}$, and

\item $\mathcal F$ contains all the diffeomorphisms
$$
f:M\rightarrow[0,a^{-1}]\times[0,a]
$$
such that
$$
f(A_1)=(0,0),\quad
f(A_2)=(a^{-1},0),\quad
f(A_3)=(a^{-1},a),\quad
f(A_4)=(0,a).
$$
\end{enumerate}
Then
\begin{enumerate}[label=(\alph*)]
\item For any $f\in\mathcal F$, $\mathcal E[f]\geq1$ and the equality holds if and only if $f$ is conformal.

\item There exists a unique harmonic map $f\in \mathcal F$, and it is the unique minimizer of the Dirichlet energy in $\mathcal F$.

\item If $f\in\mathcal F$ is conformal, then it is harmonic.
\end{enumerate}
\end{theorem}

By part (b) and (e) of Theorem 2.5, for closed surfaces it makes good sense to compute the harmonic map within a fixed homotopy class. Such a topological constrain also often naturally arises in applications.
Part (b) of Theorem 2.5 and part (b) of Theorem 2.6 somewhat justify the physical intuition that harmonic maps minimize the stretching energy. Minimizing the Dirichlet energy would also be an important idea for computing harmonic maps.
Part (c) of Theorem 2.5 and part (c) of Theorem 2.6 connect the notions of conformal maps and harmonic maps. So conformal maps are really special harmonic maps, and the converse is true if the target surface is a unit sphere. Part (a) of Theorem 2.5 and Part (a) of Theorem 2.6 point out that conformal maps minimize the normalized Dirichlet energy $\mathcal E(f)/Area(N)$ among the harmonic maps.
This optimality provides a method to compute the uniformization using harmonic maps.
\begin{enumerate}[label=(\alph*)]
\item If $M$ is a topological sphere, a harmonic diffeomorphism from $(M,g)$ to a unit sphere is already a conformal map.

\item If $M$ is a genus 1 closed orientable surface, one can first compute the harmonic map to a flat torus $(N,h)=\mathbb C/\mathbb Z\oplus\mathbb Z\tau$ for some parameter $\tau$,
and then minimize the normalized Dirichlet energy, by tuning the parameter $\tau$ in the upper half plane.

\item If $M$ is a closed orientable surface with genus $>1$, one can first compute the harmonic map to a hyperbolic surface $(N,h)=\mathbb H^2/\Gamma$ for some discrete subgroup $\Gamma$ of $Isom^+(\mathbb H^2)$,
and then minimize the normalized Dirichlet energy, by continuously deforming the subgroup $\Gamma$.

\item If $M$ is homeomorphic to a closed disk, then one can compute the harmonic map to a rectangle $[0,a^{-1}]\times[0,a]$,
and then minimize the Dirichlet energy, by tuning the parameter $a\in\mathbb R_{>0}$.
\end{enumerate}
\subsection{Discrete Harmonic Maps}
Since we will be computing discrete harmonic maps on lattice graphs, here we briefly introduce the theory of discrete harmonic maps, which has been well-studied. See \cite{hass2015simplicial}\cite{kajigaya2019uniformizing}\cite{de1991comment} for references of the related definitions and theorems.
In this subsection, we always assume that $G=(V,E)$ is an undirected connected simple graph, and $w\in\mathbb R^E_{>0}$ denotes the edge weight, and $(N,h)$ is a flat rectangle or a closed orientable surface of constant curvature $\pm1$ or 0.

A \emph{discrete map} $f$ from $G$ to $(N,h)$ maps each vertex $i\in V$ to a point in $N$, and each edge $ij\in E$ to a geodesic segment in $N$ with endpoints $f(i),f(j)$.
Given a discrete map $f$, we denote $l_{ij}(f)$ as the length of $f(ij)$. Assuming the stretching energy on a single edge $ij$ is $\frac{1}{2}w_{ij}l_{ij}^2$, we have the following definition of the discrete Dirichlet energy.
\begin{definition}
Given a discrete map $f$ from $(G,w)$ to $(N,h)$, the \emph{discrete Dirichlet energy} is defined to be
$$
\mathcal E(f)=\frac{1}{2}\sum_{ij\in E}w_{ij}\cdot l_{ij}(f)^2.
$$
\end{definition}

%
%
%
%

The notion of discrete harmonic maps can be defined by the balanced conditions on vertices.
\begin{definition}
A discrete map $f$ is \emph{harmonic} if for any $i\in V$
\begin{equation}
\sum_{j:ij\in E}w_{ij}\cdot l_{ij}(f)\cdot\vec{e}_{ij}=0,
\end{equation}
where $\vec e_{ij}$ is the unit vector in $T_{f(i)}N$ representing the direction of the geodesic $f(ij)$.
\end{definition}
This definition consists with the continuous one in the following sense.
\begin{proposition}A discrete harmonic map $f$ is a critical point of $\mathcal E$, if $N$ is not a sphere or $l_{ij}(f)<\pi$ for any $ij\in E$.
\end{proposition}
Here the assumption on $N$ and $l_{ij}(f)$ is to guarantee that locally $f$ can be parameterized by $[f(i)]_{i\in V}\subset N^V$.
For the cases of non-positive curvatures, we have
the following two well-known theorems partially analogue to Theorem 2.5 and 2.6.

\begin{theorem}If $(N,h)$ is closed and has constant curvature $0$ or $-1$, then
any discrete map $f_0$ from $G$ to $N$ is homotopic to a discrete harmonic map $f$, which is a minimizer of the discrete Dirichlet energy $\mathcal E$ in the homotopy class.
Such discrete harmonic map $f$ is

\begin{enumerate}[label=(\alph*)]
\item unique up to a translation, if $(N,h)$ is a flat torus, and

\item unique if $(N,h)$ has constant curvature $-1$.
\end{enumerate}

\end{theorem}

\begin{theorem}Assume $(N,h)$ is a flat rectangle with four ordered edges $B_1,B_2,B_3,B_4$, and $V_1,V_2,V_3,V_4$ are 4 subsets of $V$ such that $V_1\cap V_3=\emptyset$ and $V_2\cap V_4=\emptyset$. Then there exists a unique discrete map $f$ from $G$ to $N$ such that
\begin{equation*}
f(V_1)\subset B_1,\quad
f(V_2)\subset B_2,\quad
f(V_3)\subset B_3,\quad
f(V_4)\subset B_4,
\end{equation*}
and the balanced condition (1) holds for any $i\in V-V_1\cup V_2\cup V_3\cup V_4$. Such constrained discrete harmonic map minimizes the discrete Dirichlet energy $\mathcal E(f)$ among all the constrained maps.
\end{theorem}

\section{Cubic Lattice Approximation}
\subsection{Construction of Lattices}

Assume that $M$ is a closed smooth surface embedded in $\mathbb R^3$, and $P$ is a $\delta$-point cloud of $M$. If $\delta$ is sufficiently small, comparing to some constant $\epsilon>0$, then the $\epsilon$-neighborhood
$B(P,\epsilon)$ of $P$ would be a good approximation of the $\epsilon$-neighborhood $B(M,\epsilon)$ of $M$. If $n$ is sufficiently large, comparing to $1/\epsilon$, then
$$
V=B(P,\epsilon)\cap(\mathbb Z/n)^3
$$
would be a good cubic lattice approximation of $B(P,\epsilon)\approx B(M,\epsilon)$ (see Figure 1). Connect
two vertices $x,y$ in $V$ by an edge if $d(x,y)=1/n$, and then
we obtain a lattice graph $G=(V,E)$.
Since we are using standard cubic lattice, the weight function is set to be constant $w_{ij}=1$. Such weighted graph $(G,w)$ would be our lattice approximation of $M$.

In practice, if the pointcloud data $P$ is given, we need to choose a proper $\epsilon>0$ such that
\begin{enumerate}[label=(\alph*)]
\item $\epsilon$ is sufficiently large such that $B(P,\epsilon)$ forms a connected domain in $\mathbb R^3$ and has uniform \emph{thickness}, and

\item $\epsilon$ is sufficiently small such that $B(P,\epsilon)\approx B(M,\epsilon)$ is a good approximation of $M$.
\end{enumerate}
For our second parameter $n$, if we choose a larger integer $n$, the lattice graph $G$ should approximates the domain $B(P,\epsilon)\approx B(M,\epsilon)$ better. However, in practice we find that sparse lattices already work well.

\begin{figure}[H]
\centering
\includegraphics[width=0.25\textwidth]{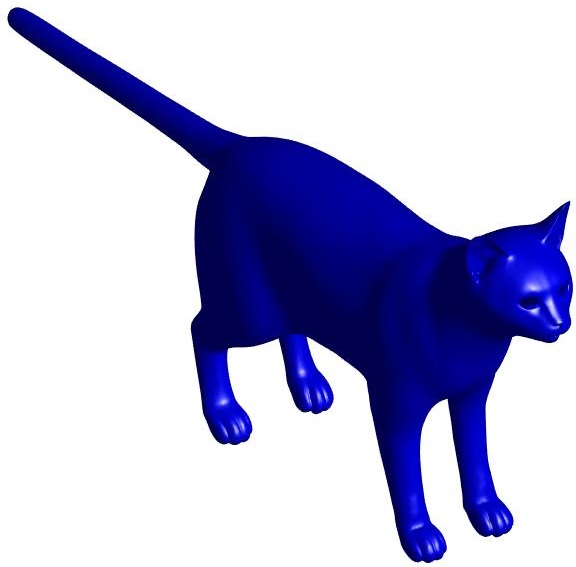}\quad\quad\quad
\includegraphics[width=0.25\textwidth]{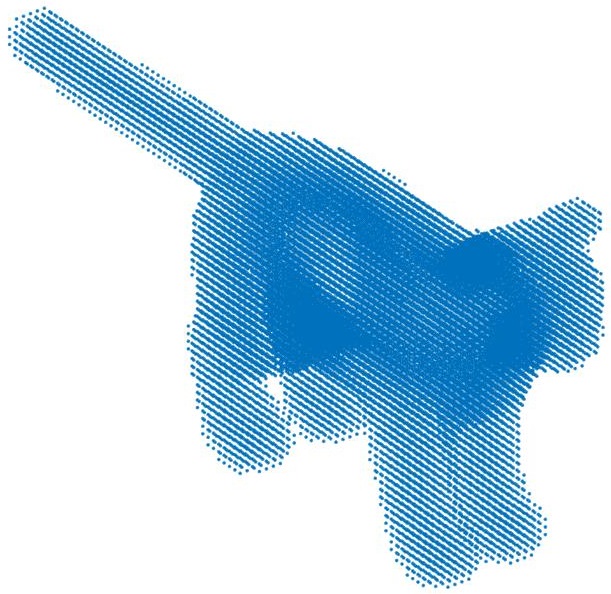}
\caption{Lattice approximation of a cat}
\end{figure}

\subsection{Trilinear interpolation}
Suppose $G=(V,E)$ is a lattice approximation of $P$ constructed as before, then given a function $f:V\rightarrow \mathbb R^k$, we can use the standard trilinear interpolation to extend $f$ to the point cloud $P$. Assume $p=(p_1,p_2,p_3)\in P$ lies inside a cubic cell $[x_0,x_1]\times[y_0,y_1]\times[z_0,z_1]$ of the lattice $G$. Denote $a_{ijk}=f(x_i,y_j,z_k)$ where $i,j,k=\in\{0,1\}$, and
\begin{align*}
x_0'=p_1-x_0,\quad
y_0'&=p_2-y_0,\quad
z_0'=p_3-z_0,\\
x_1'=x_1-p_1,\quad
y_1'&=y_1-p_2,\quad
z_1'=z_1-p_3,\\
\end{align*}
then
$$
f(p)=n^3\sum_{i,j,k\in\{0,1\}}a_{ijk}x'_iy'_jz'_k
$$
is the trilinear interpolation of $f$ on $P$.

\section{Computing Harmonic and Conformal Maps to Spheres}

By part (c)(d) of Theorem 2.5, in the spherical case harmonic maps and conformal maps are really the same. This section solves the following problem of computing harmonic and conformal maps.

\begin{problem}
Assume that $M$ is a genus 0 closed smooth surface embedded in $\mathbb R^3$.
Given a point cloud approximation $P$ of $M$, how to compute harmonic diffeomorphisms, i.e., conformal diffeomorphisms from $M$ to the unit sphere?
\end{problem}

\subsection{Algorithm}
Our method is smoothly adapted from Gu-Yau's method \cite{gu2002computing} for computing harmonic maps to spheres.
First we pick a lattice approximation $G=(V,E)$ as in Section 3, and then
\begin{enumerate}[label=(\arabic*)]
\item compute the (approximated) discrete harmonic map $f:G\rightarrow\mathbb S^2$, and

\item extend $f$ to $P$ by trilinear interpolation, and then do a normalization $x\mapsto x/|x|_2$ and get
an approximation of a harmonic map from $M$ to the unit sphere.
\end{enumerate}

Step (2) is pretty clear and simple, so let us discuss the details of step (1).
Recall that a discrete map $f$ from $G$ to $\mathbb S^2\subset \mathbb R^3$ is harmonic if it is a critical point of the energy
$$
\mathcal E(f)=\frac{1}{2}\sum_{ij}w_{ij}l_{ij}^2(f).
$$
If the edge length $l_{ij}(f)$ is small, $l_{ij}(f)\approx |f(j)-f(i)|_2$, and the discrete energy is approximated by
$$
\mathcal E_0(f)=\frac{1}{2}\sum_{ij\in E}w_{ij}|f(j)-f(i)|_2^2.
$$
So for simplicity we will minimize $\mathcal E_0(f)$ instead of $\mathcal E(f)$ to compute the approximation of the discrete harmonic map.
Since
$$
\frac{\partial \mathcal E_0}{\partial f(i)}=-(Lf)(i),
$$
a critical point $f$ of $\mathcal E_0$ satisfies
that for any $i\in V$
$$
(Lf)(i)\perp T_{f(i)}\mathbb S^2,
$$
i.e.,
\begin{equation}
(L^{\|}f)(i)=0
\end{equation}
where
$
L^{\|}f(i)
$
denotes
the tangential component of $Lf(i)$, i.e.,
the orthogonal projection of $Lf(i)$ to the plane $f(i)^{\perp}$.
We view a map $f$ satisfying equation (2) as an approximation of a discrete harmonic map, and wish to compute it by the following discrete heat flow on $\mathbb S^2$
$$
\frac{df}{dt}=-L^\|f.
$$
We can simply solve the above ODE by the explicit Euler's method, with a normalization after each iteration to keep that all the points are on the unit sphere. More specifically, if the step size is $\delta t$ we have that
\begin{equation*}
f^{t+\delta t}=\pi\circ(f^t-\delta t\cdot L^\| f^t)
\end{equation*}
where $\pi(x)=x/|x|_2$.

Two harmonic maps to a sphere can differ by a M\"obius transformation.
So a smooth harmonic map to a sphere is not unique, and would
have large distortion if most of the mass is concentrated near a point of the sphere. To avoid such big distortion, we add a correction mapping $m$ in each iteration. Here the correction mapping $m$ first translate all the points $f(i)$'s in $\mathbb R^3$ to make the center of mass be at the origin, and then do a normalization $x\mapsto x/|x|_2$.
Now the modified iteration is
\begin{equation}
f^{t+\delta t}=m\circ \pi\circ(f^t-\delta t\cdot L^\| f^t),
\end{equation}
and we have the following algorithm.

\begin{algorithm}
\caption{Harmonic and Conformal Maps to Spheres}

\textbf{Input:}  A point cloud $P$ in $\mathbb R^3$\\
\textbf{Output:}  A harmonic map $f$ from $P$ to the unit sphere

\begin{algorithmic}[1]
\State Translate $P$ such that the origin lies inside of $P$.
\State Choose parameters $\epsilon$ and $n$, and then compute the lattice approximation $G=(V,E)$.
\State Compute $f(x)=\pi(x)$ on $V$ as the initial map.
\State Do the iteration by equation (3), until the function $f$ converges.
\State Extend $f$ to $P$ by trilinear interpolation.
\State \textbf{Return} $\pi\circ (f|_P)$.

\end{algorithmic}
\end{algorithm}

\subsection{Numerical Examples}
See the two numerical examples in Figure 2 and Figure 3.

\begin{figure}[H]
\centering
\includegraphics[width=0.35\textwidth]{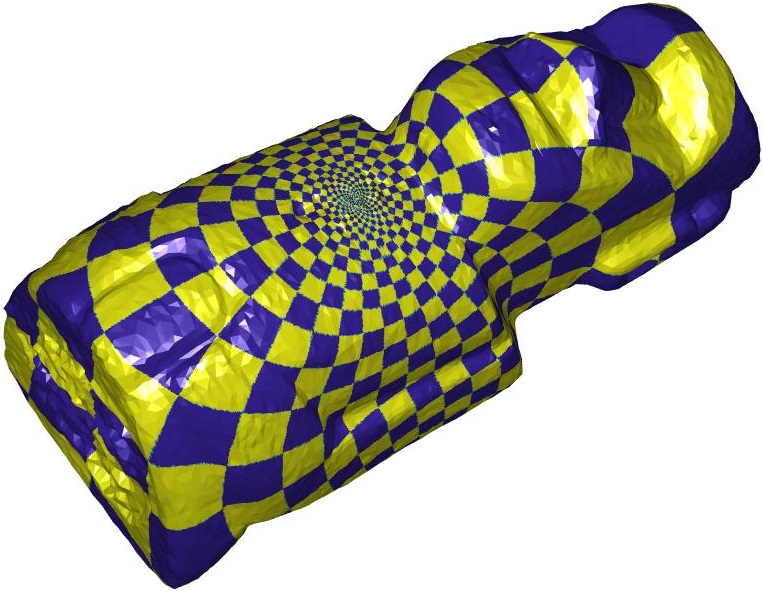}\quad
\includegraphics[width=0.25\textwidth]{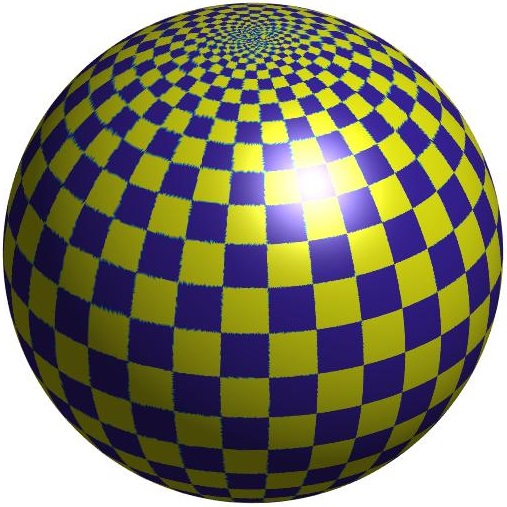}
\caption{Harmonic map from a moai to a unit sphere}
\end{figure}
\begin{figure}[H]
\centering
\includegraphics[width=0.35\textwidth]{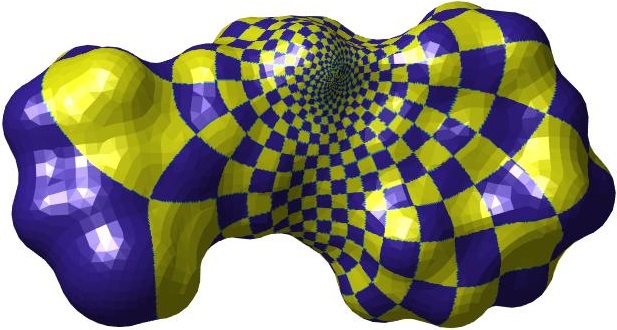}\quad
\includegraphics[width=0.25\textwidth]{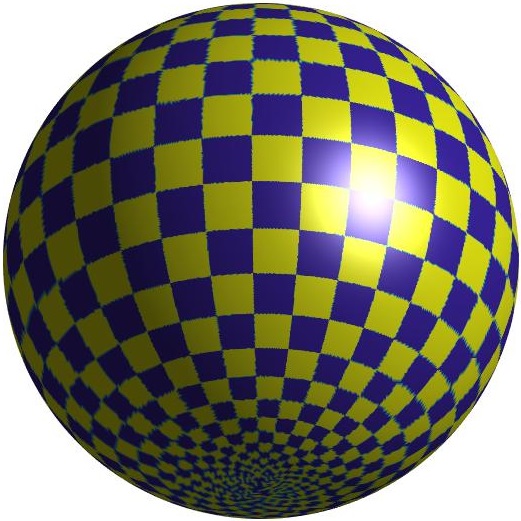}
\caption{Harmonic map from retinal to a unit sphere}
\end{figure}

\section{Computing Harmonic and Conformal Maps to Rectangles}
This section solves the following problem on computing harmonic and conformal maps.
\begin{problem}
Assume that $M$ is a smooth surface embedded in $\mathbb R^3$, and is homeomorphic to a closed disk, and
$\gamma_1,\gamma_2,\gamma_3,\gamma_4$ are 4 ordered arcs on $\partial M$ divided by 4 points on $\partial M$.
Now we are given the data $(P,P_1.P_2,P_3,P_4)$ such that
$P$ is a point cloud approximating $M$, and $P_i$ is a subset of $P$ approximating $\gamma_i$.
Here are two computational problems.
\begin{enumerate}[label=(\arabic*)]
\item Compute the harmonic map $f:M\rightarrow[0,a^{-1}]\times[0,a]$ for a given $a>0$ such that
\begin{align}
f(\gamma_1)=&[0,a^{-1}]\times\{0\},\\
f(\gamma_2)=&\{a^{-1}\}\times[0,a],\\
f(\gamma_3)=&[0,a^{-1}]\times\{a\},\\
f(\gamma_4)=&\{0\}\times[0,a].
\end{align}

\item Find a positive number $a>0$ and a conformal map $f:M\rightarrow[0,a^{-1}]\times[0,a]$ such that the above boundary conditions (4)(5)(6)(7) are satisfied.
\end{enumerate}

\end{problem}
We will first discuss the harmonic maps in Section 5.1, and then discuss the conformal maps in Section 5.2.
\subsection{Algorithm for Harmonic Maps}
First we construct a lattice approximation $G=(V,E)$ of $P$ as in Section 3, and let $V_i=\cup_{x\in P_i}V_x$ be the lattice approximation of $\gamma_i$, where $V_x$ consists of the 8 vertices of the cubic cell of $G$ containing $x\in P$.
Then we would like to
\begin{enumerate}[label=(\arabic*)]
\item compute the (constrained) discrete harmonic map $f:G\rightarrow\mathbb [0,a^{-1}]\times[0,a]$, such that
\begin{align*}
f(V_1)\subset&[0,a^{-1}]\times\{0\},\\
f(V_2)\subset&\{a^{-1}\}\times[0,a],\\
f(V_3)\subset&[0,a^{-1}]\times\{a\},\\
f(V_4)\subset&\{0\}\times[0,a],
\end{align*}
and

\item extend $f$ to $P$ by trilinear interpolation, and then $f|_P$ is our approximation of the desired  harmonic map.
\end{enumerate}

Step (2) is clear and simple. Computing $f=(f^1,f^2)$ in step (1) amounts to solve the following 2 systems of linear equations
\begin{equation}
\left\{
\begin{array}{lll}
f^1(i)=0&\quad\text{ if }& i\in V_4\\
f^1(i)=a^{-1}&\quad\text{ if }& i\in V_2\\
Lf^1(i)=0&\quad\text{ if }& i\in V-V_4-V_2
\end{array}
\right.
,
\left\{
\begin{array}{lll}
f^2(i)=0&\quad\text{ if }& i\in V_1\\
f^2(i)=a&\quad\text{ if }& i\in V_3\\
Lf^2(i)=0&\quad\text{ if }&i\in V-V_1-V_3
\end{array}
\right..
\end{equation}
These two equations can be solved efficiently by the preconditioned conjugate gradient method.

In a summary we have the following \textbf{Algorithm 2}.

\begin{algorithm}
\caption{Harmonic Maps to Rectangles}

\textbf{Input:}  A point cloud $P$ in $\mathbb R^3$ as an approximation of a surface $M$, and 4 subsets $P_1,P_2,P_3,P_4$ as approximations of 4 boundary arcs of $M$, and a parameter $a\in\mathbb R_{>0}$.
\\
\textbf{Output:}  A harmonic map $f$ from $P$ to the rectangle $[0,a^{-1}]\times[0,a]$.

\begin{algorithmic}[1]
\State  Choose parameters $\epsilon$ and $n$, and then compute the lattice approximation $(G,V_1,V_2,V_3,V_4)$ of $(P,P_1,P_2,P_3,P_4)$.

\State  Solve the two systems of linear equations (8).

\State Extend $f$ to $P$ by trilinear interpolation

\State \textbf{Return} $f|_P$.

\end{algorithmic}
\end{algorithm}

\subsection{Algorithm for Conformal Maps}
Assume $f_a$ is the discrete harmonic map from the lattice approximation $G=(V,E)$ to the rectangle $[0,a^{-1}]\times[0,a]$.
To compute conformal maps, we only need to find a proper edge length $a$ such that the computed discrete harmonic map $f_a$ approximates a conformal mapping. Inspired by part (a) of Theorem 2.6, we compute the conformal map by minimizing $\mathcal E(f_a)$ over $a\in\mathbb R_{>0}$. If $\bar a$ is the minimizer, then
$
f_{\bar a}:V \rightarrow[0,\bar a^{-1}]\times[0,\bar a]
$
would approximate a conformal map.

Our method is summarized as \textbf{Algorithm 3}.

\begin{algorithm}
\caption{Conformal Maps to Rectangles}

\textbf{Input:}  A point cloud $P$ in $\mathbb R^3$ as an approximation of a surface $M$, and 4 subsets $P_1,P_2,P_3,P_4$ as approximations of 4 boundary arcs of $M$.
\\
\textbf{Output:} A parameter $a$, and a conformal map $f$ from $P$ to the rectangle $[0,a^{-1}]\times[0,a]$.

\begin{algorithmic}[1]
\State  Choose parameters $\epsilon$ and $n$, and then compute the lattice approximation $(G,V_1,V_2,V_3,V_4)$ of $(P,P_1,P_2,P_3,P_4)$.

\State  For a given $a\in\mathbb R_{>0}$, solve the two systems of linear equations (8) and get a discrete harmonic map $f_a$.

\State Compute the discrete Dirichlet energy $\mathcal E(f_a)$.

\State Compute the minimizer $\bar a$ of $\mathcal E(f_a)$, and $f_{\bar a}$.

\State Extend $f_{\bar a}$ to $P$ by trilinear interpolation.

\State \textbf{Return} $\bar a$ and $f_{\bar a}|_P$.

\end{algorithmic}
\end{algorithm}

\subsection{Numerical Examples}
 Here are two numerical examples for harmonic maps and conformal maps respectively.
\begin{figure}[H]
\centering
\includegraphics[width=0.15\textwidth]{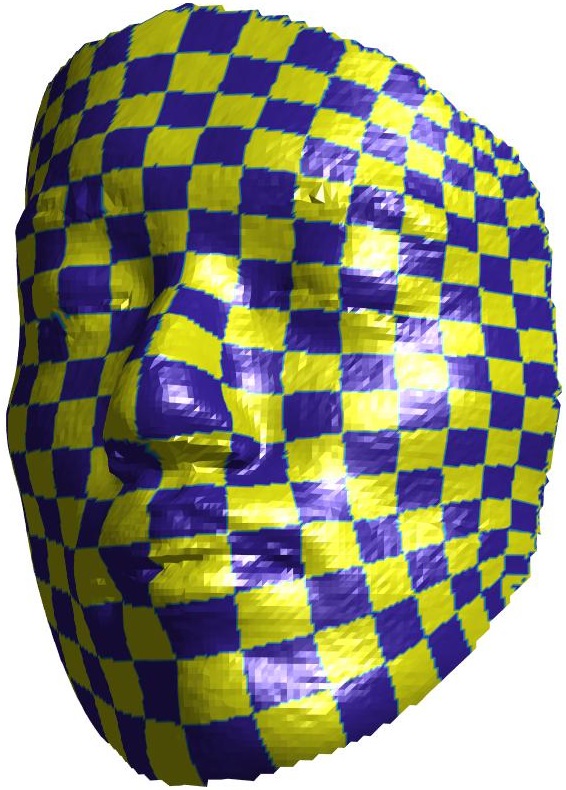}
\includegraphics[width=0.2\textwidth]{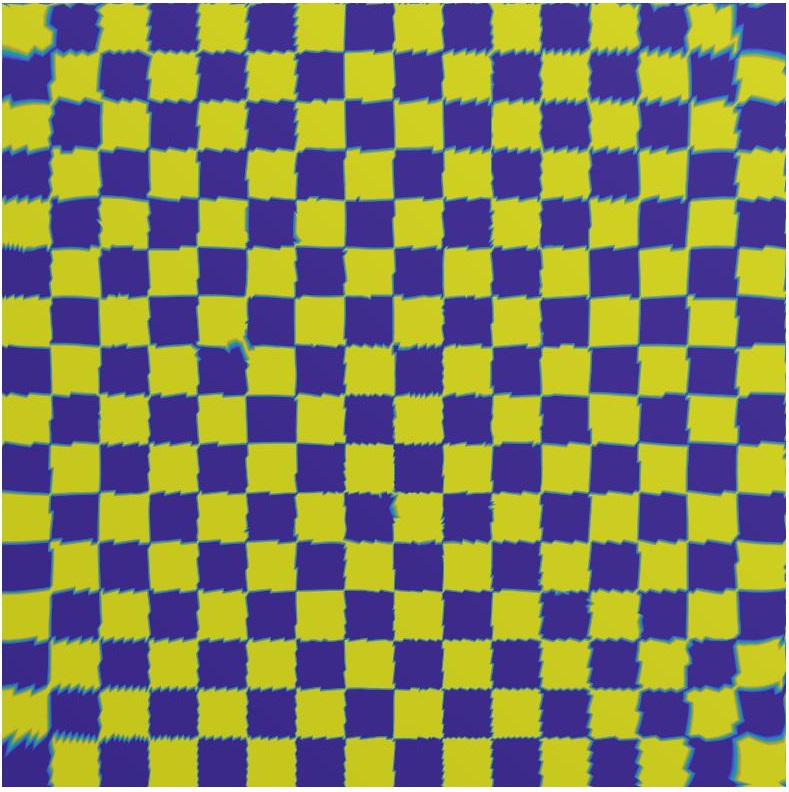}\quad\quad
\includegraphics[width=0.3\textwidth]{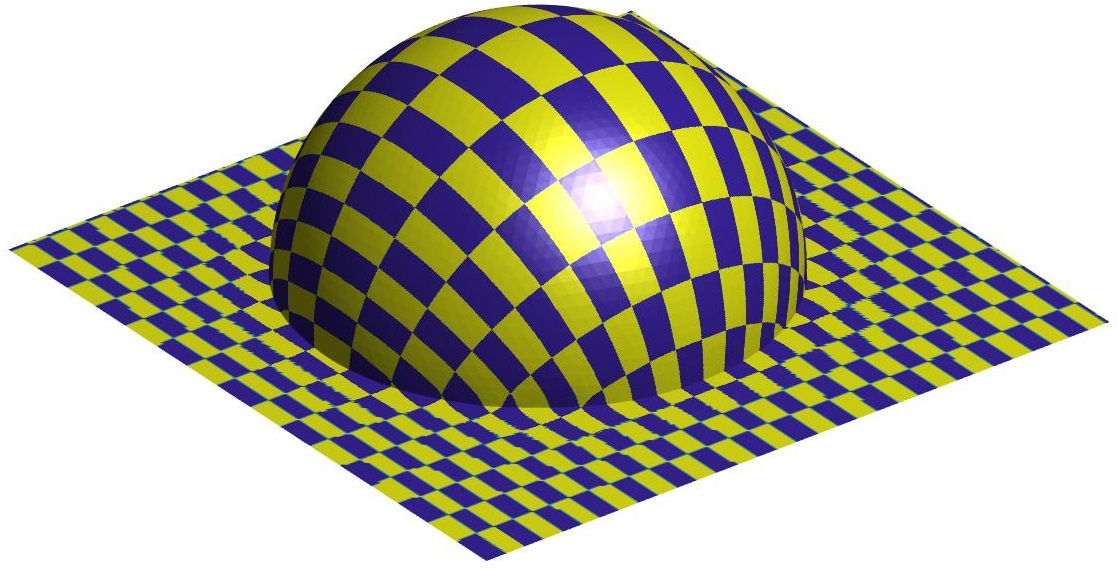}
\includegraphics[width=0.25\textwidth]{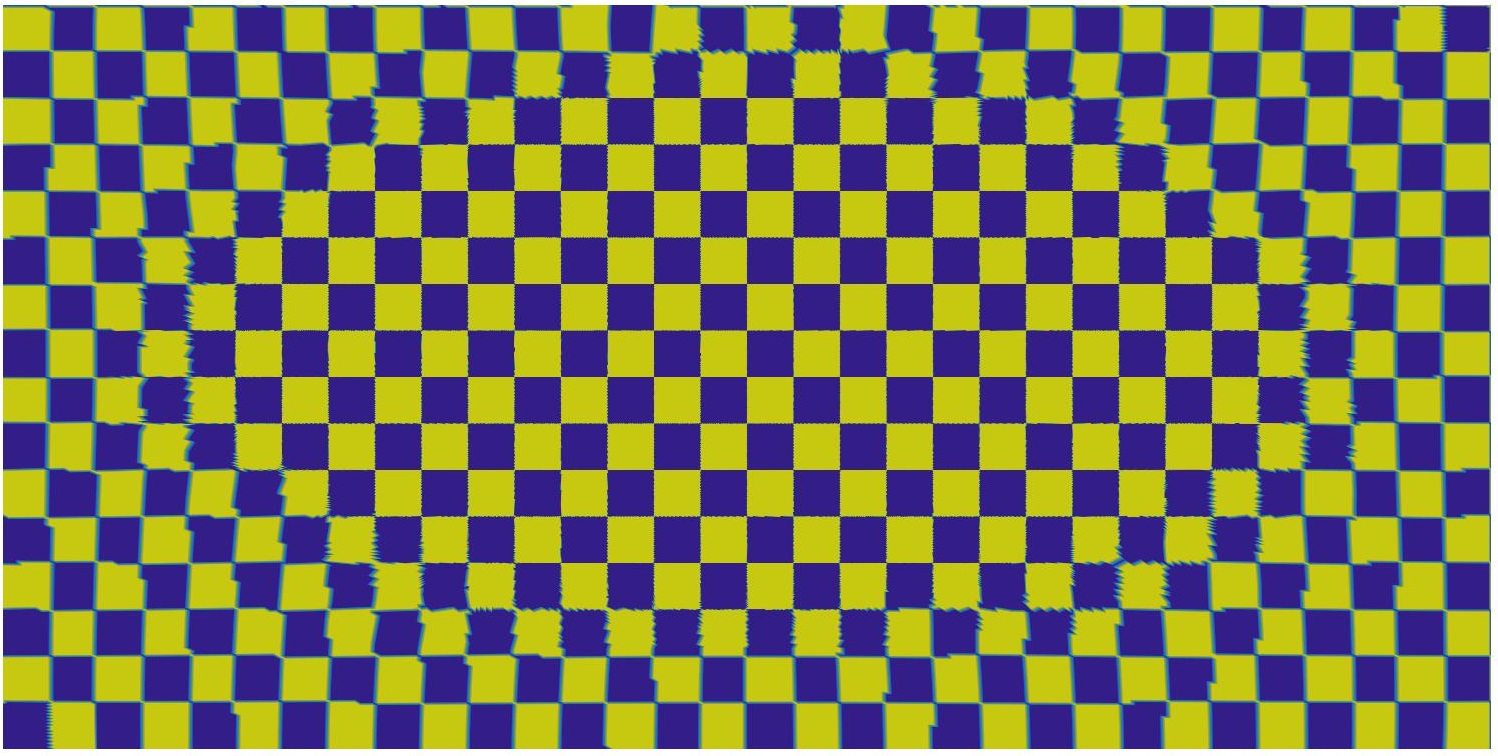}
\caption{Harmonic maps to rectangles}
\end{figure}

\begin{figure}[H]
\centering
\includegraphics[width=0.15\textwidth]{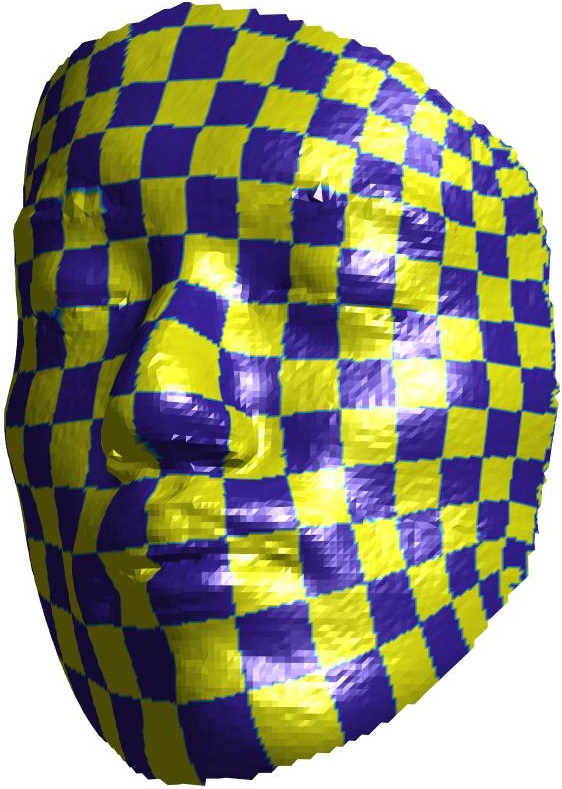}
\includegraphics[width=0.2\textwidth]{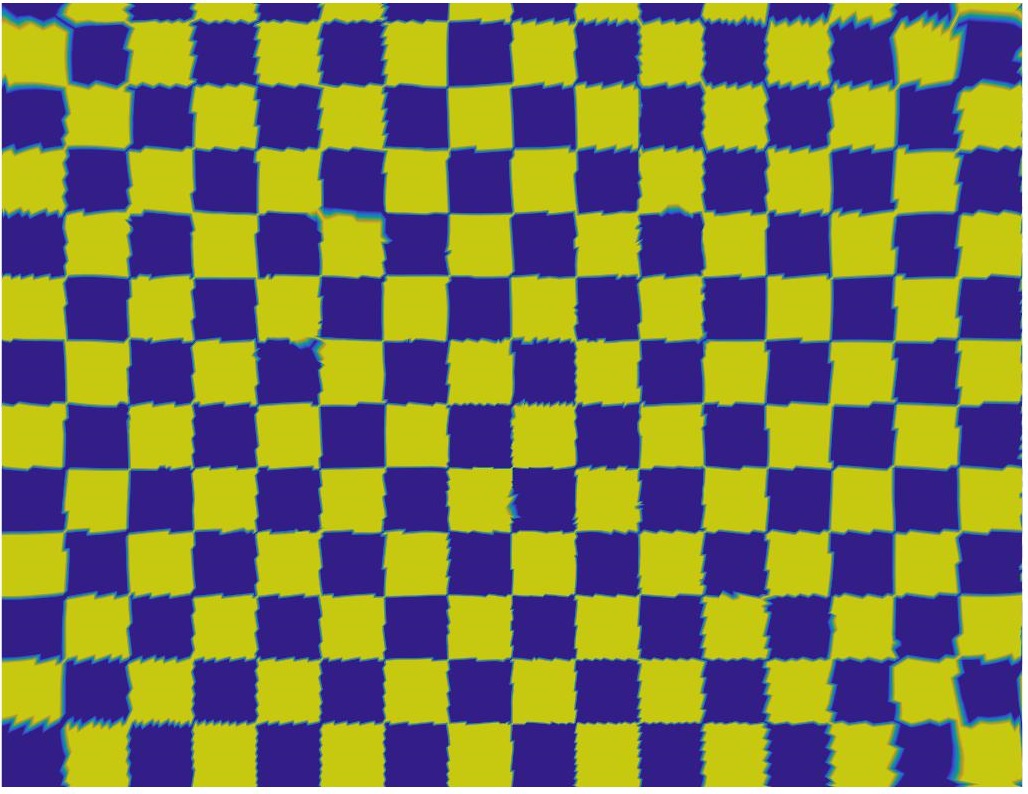}
\quad\quad
\includegraphics[width=0.3\textwidth]{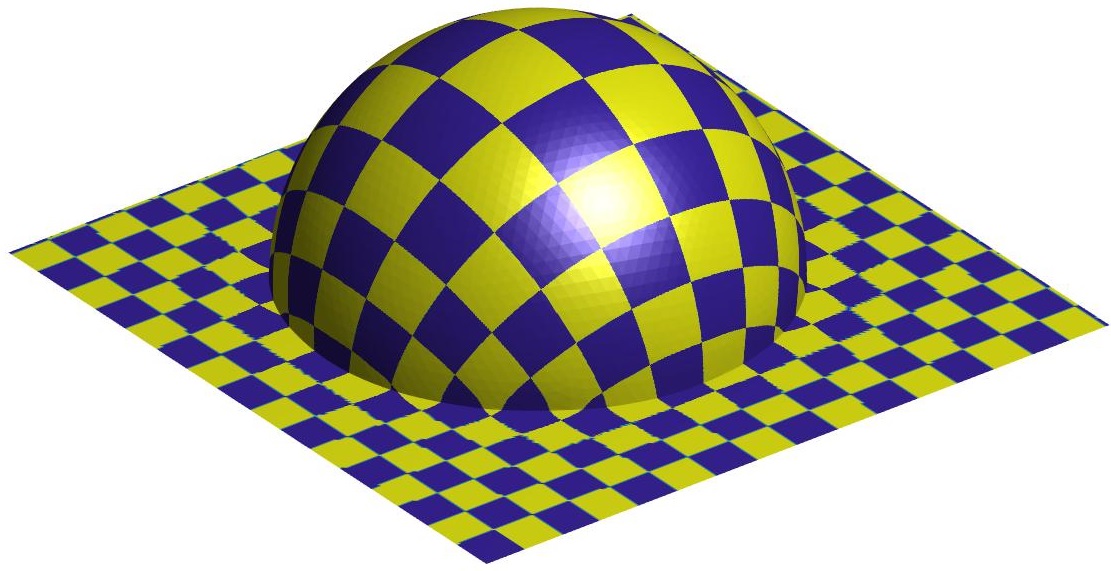}
\includegraphics[width=0.2\textwidth]{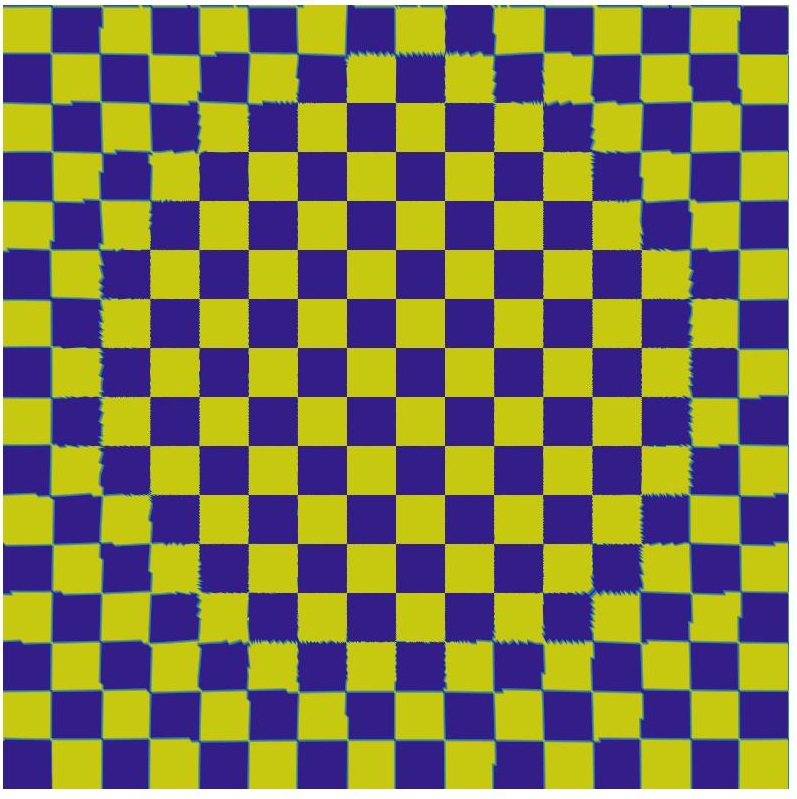}
\caption{Conformal maps to rectangles}
\end{figure}

\section{Computing Harmonic and Conformal Maps to Flat Tori}
This section solves the following problem on computing harmonic and conformal maps.
\begin{problem}
Assume that $M$ is a genus 1 closed smooth surface embedded in $\mathbb R^3$, and we are given
a point cloud approximation $P$ of $M$.
\begin{enumerate}[label=(\arabic*)]
\item How to compute a harmonic map
$$
f:M\rightarrow\mathbb C/\mathbb Z\oplus\mathbb Z\tau
$$
for a given $\tau$ with $Im(\tau)>0$?

\item How to find the complex number $\tau$ with $Im(\tau)>0$ and the map
$$
f:M\rightarrow\mathbb C/\mathbb Z\oplus\mathbb Z\tau
$$
such that $f$ is conformal?
\end{enumerate}
\end{problem}
We will first discuss the harmonic maps in Section 6.1, and then discuss the conformal maps in Section 6.2.
\subsection{Algorithm for Harmonic Maps}

First we construct a lattice approximation $G=(V,E)$ as in Section 3, and then need to

(1) compute the discrete harmonic map $f$ from $G$ to $\mathbb C/\mathbb Z\oplus\mathbb Z\tau$, or equivalently, to a fundamental domain in $\mathbb C$, and

(2) extend $f$ to $P$ by trilinear interpolation, and then $f|_P$ represents
an approximation of the desired harmonic map.

Step (2) is kind of clear and simple, so let us discuss the details of step (1).
First we "cut" the lattice $G$ along two "loops" $\gamma_1,\gamma_2$, so that we can represent a discrete map from $G$ to a torus as a map from $V$ to a planar domain in $\mathbb C$. Secondly we fix the homotopy class of discrete maps $f:G\rightarrow\mathbb C/\mathbb Z\oplus\mathbb Z\tau$, by requiring that $\gamma_1$ corresponds to the translation $z\mapsto z+1$, and $\gamma_2$ corresponds to the translation $z\mapsto z+\tau$.
Now computing the discrete harmonic map in the fixed homotopy class amounts to solve a system of complex linear equations. For a vertex $i$ that is away from the loops $\gamma_1$ and $\gamma_2$, the balanced condition (1) gives us that
$$
Lf(i)=0.
$$
Things become a bit subtle when the vertex $i$ is near a loop. Assume the edge $ij$ pass through the loop $\gamma_1$, and $l_{ij}(f)$ is small, and all the other edges adjacent to vertex $i$ do not pass through $\gamma_1$ or $\gamma_2$. Then
$$
f(j)\approx f(i)\pm\tau,
$$
where the sign depends on the orientation of $\gamma_1$ and the relative position between $\gamma_1$ and $i$. For such a vertex $i$, the corrected balanced condition should be
$$
Lf(i)=\pm\tau.
$$
Similar corrections should be made for all the edges passing through $\gamma_1$ or $\gamma_2$. After making all the corrections, we arrived at a system of complex linear equations of the form
\begin{equation}
Lf(i)=b(i),\quad\forall i
\end{equation}
where $b(i)$ is computed by adding up all the correction terms. The equation (9) can be solved efficiently by the preconditioned conjugate gradient method.

In a summary we have the following \textbf{Algorithm 4}.

\begin{algorithm}
\caption{Harmonic Maps to Flat Tori}

\textbf{Input:}  A point cloud $P$ in $\mathbb R^3$ as an approximation of a surface $M$, and a parameter $\tau\in\mathbb C$ with $Im(\tau)>0$.
\\
\textbf{Output:}  A harmonic map $f$ from $P$ to the flat torus $\mathbb C/\mathbb Z\oplus\mathbb Z\tau$.

\begin{algorithmic}[1]
\State  Choose parameters $\epsilon$ and $n$, and then compute the lattice approximation $G=(V,E)$.

\State Choose two "loops" $\gamma_1,\gamma_2$ in $G$.

\State Compute the total correction term $b(i)$'s.

\State Solve the linear systems (9).

\State Extend $f$ to $P$ by trilinear interpolation.

\State \textbf{Return} $f|_P$.

\end{algorithmic}
\end{algorithm}
\subsection{Algorithm for Conformal Maps}
Assume that $f_\tau$ is the discrete harmonic map from the lattice approximation $G=(V,E)$ to the flat torus $\mathbb C/\mathbb Z\oplus\mathbb Z\tau$.
Inspired by part (a) of Theorem 2.5, we compute the conformal map by minimizing
$$
\frac{\mathcal E(f_\tau)}{Area(\mathbb C/\mathbb Z\oplus\mathbb Z\bar\tau)}=\frac{\mathcal E(f_\tau)}{Im(\tau)}
$$
over the parameter $\tau$. If $\bar\tau$ is the minimizer, then
$
f_{\bar\tau}:V \rightarrow\mathbb C/\mathbb Z\oplus\mathbb Z{\bar\tau}
$
would approximate a conformal map.

Our method is summarized as \textbf{Algorithm 5}.

\begin{algorithm}
\caption{Conformal Maps to Flat Tori}

\textbf{Input:}  A point cloud $P$ in $\mathbb R^3$ as an approximation of a surface $M$.
\\
\textbf{Output:} A parameter $\tau\in\mathbb C$ with $Im(\tau)>0$, and a conformal map $f$ from $P$ to the flat torus $\mathbb C/\mathbb Z\oplus\mathbb Z\tau$.

\begin{algorithmic}[1]
\State  Choose parameters $\epsilon$ and $n$, and then compute the lattice approximation $G=(V,E)$.

\State Choose two "loops" $\gamma_1,\gamma_2$ in $G$.

\State Compute the total correction term $b(i)$'s.

\State Solve the linear systems (9).

\State Compute the normalized discrete Dirichlet energy $\mathcal E(f_\tau)/Im(\tau)$.

\State Compute the minimizer $\bar \tau$ of $\mathcal E(f_\tau)/Im(\tau)$, and $f_{\bar \tau}$.

\State Extend $f_{\bar \tau}$ to $P$ by trilinear interpolation.

\State \textbf{Return} $\bar \tau$ and $f_{\bar \tau}|_P$.
\end{algorithmic}
\end{algorithm}

\subsection{Numerical Examples}
Here are two numerical examples for harmonic maps and conformal maps respectively.
\begin{figure}[H]
\centering
\includegraphics[width=0.15\textwidth]{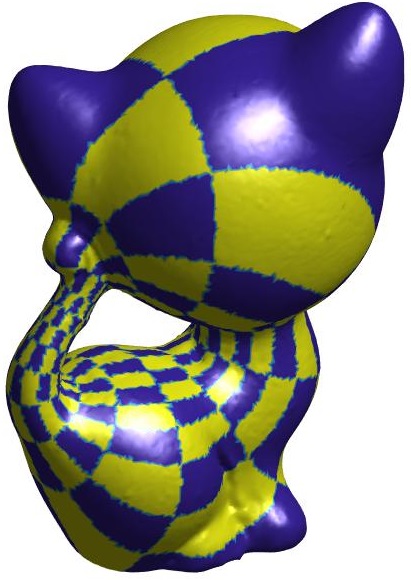}
\includegraphics[width=0.2\textwidth]{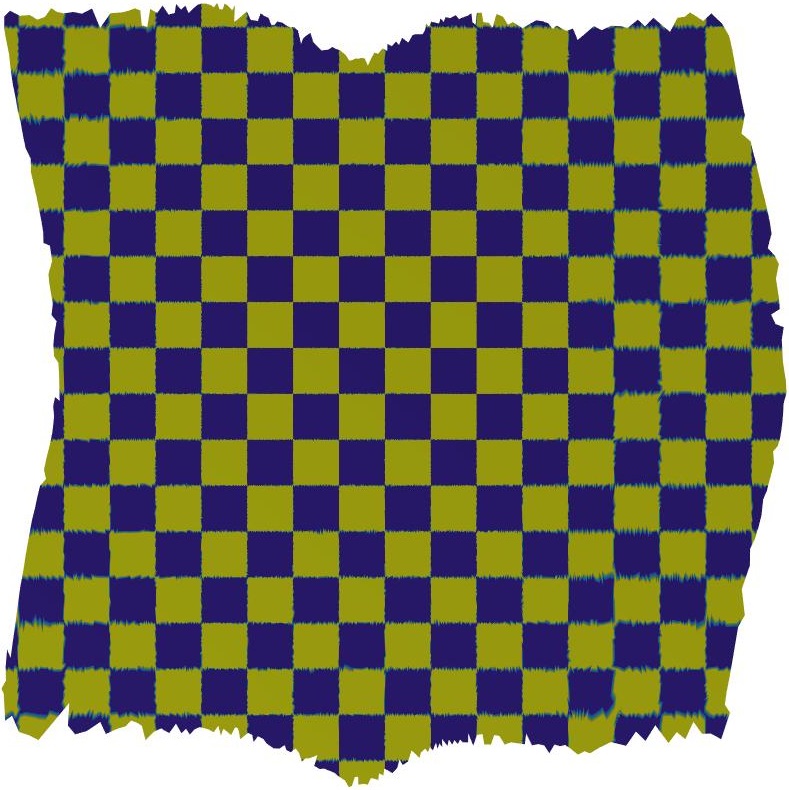}
\quad\quad
\includegraphics[width=0.2\textwidth]{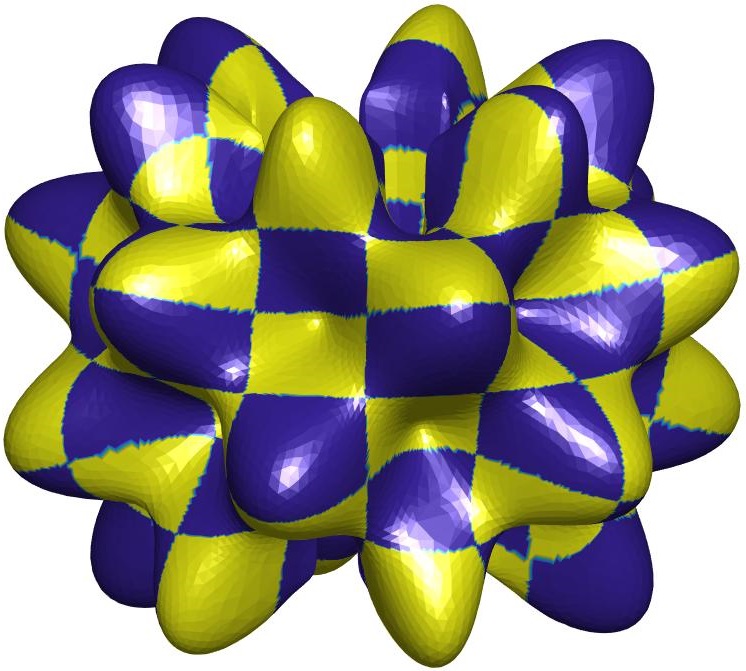}
\includegraphics[width=0.2\textwidth]{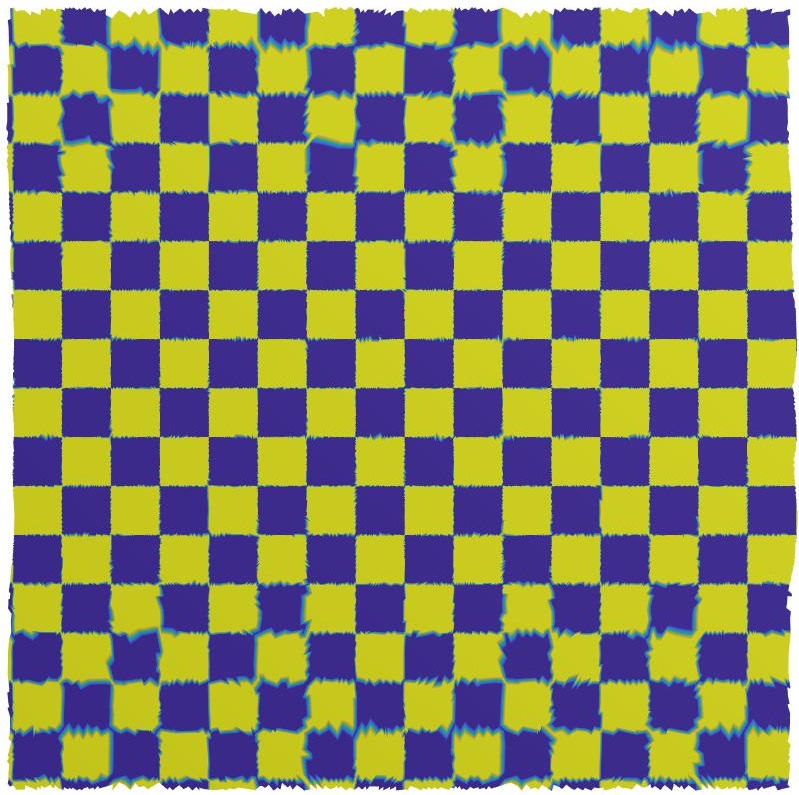}
\caption{Harmonic maps to flat tori}
\end{figure}

\begin{figure}[H]
\centering
\includegraphics[width=0.15\textwidth]{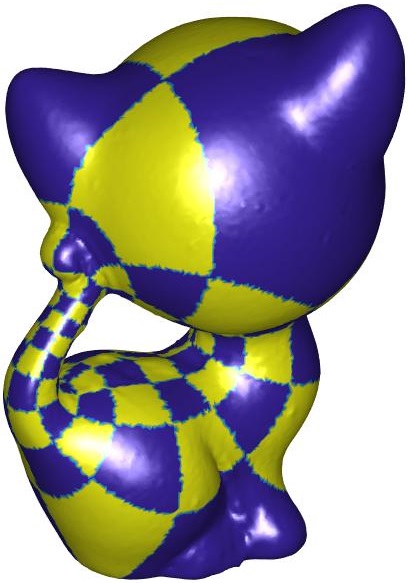}
\includegraphics[width=0.25\textwidth]{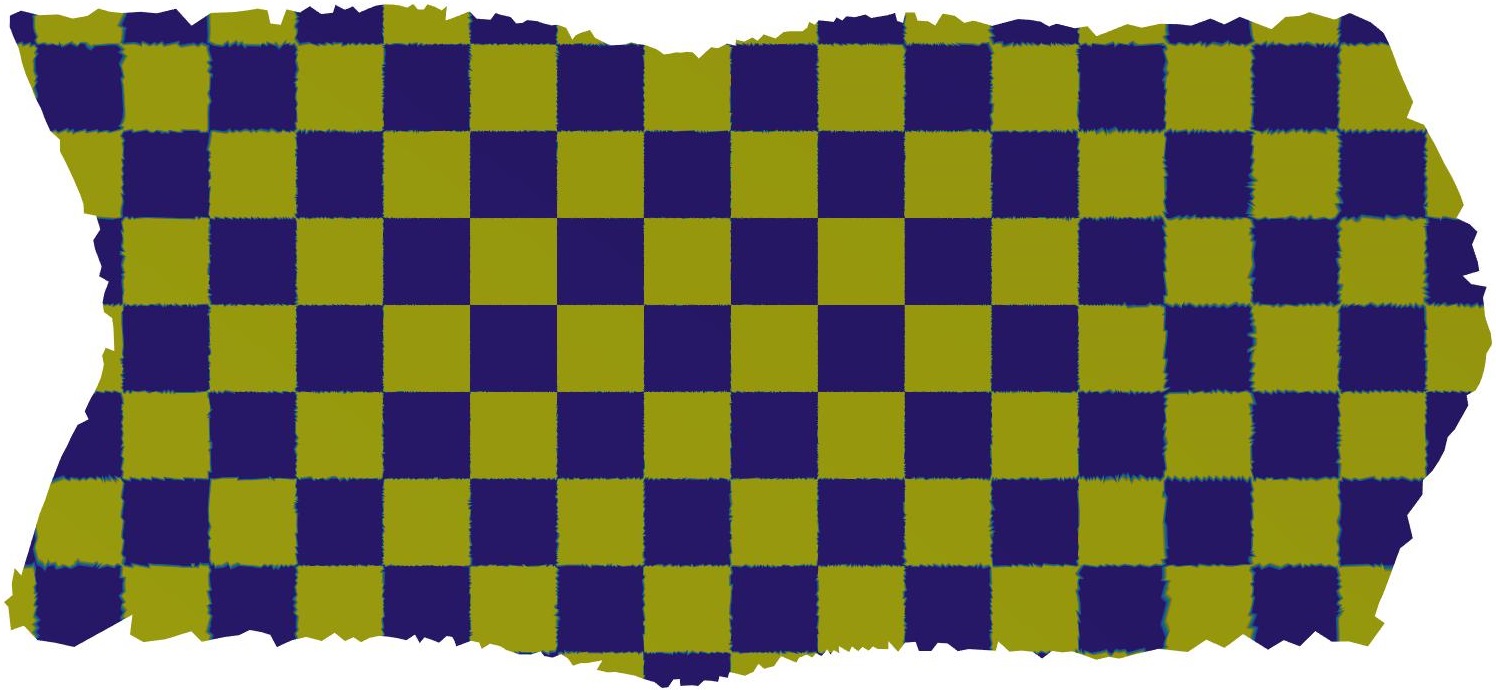}
\quad\quad
\includegraphics[width=0.2\textwidth]{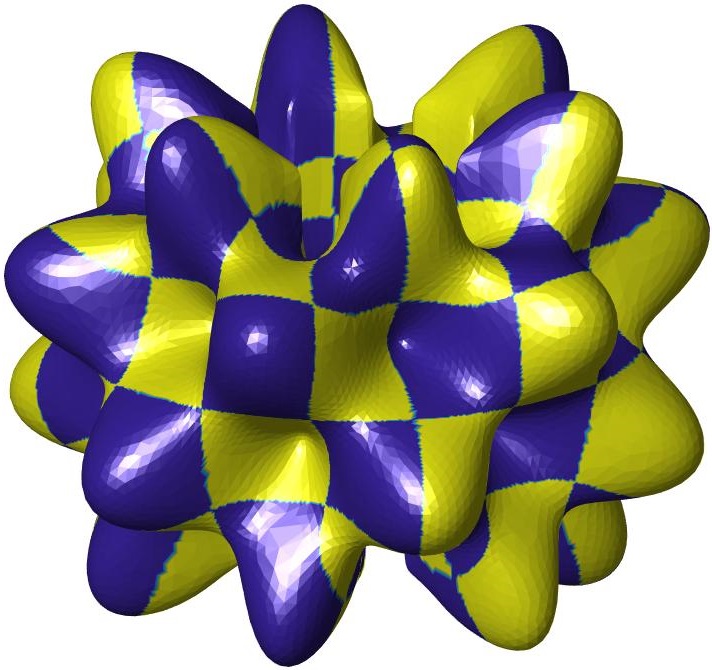}
\includegraphics[width=0.2\textwidth]{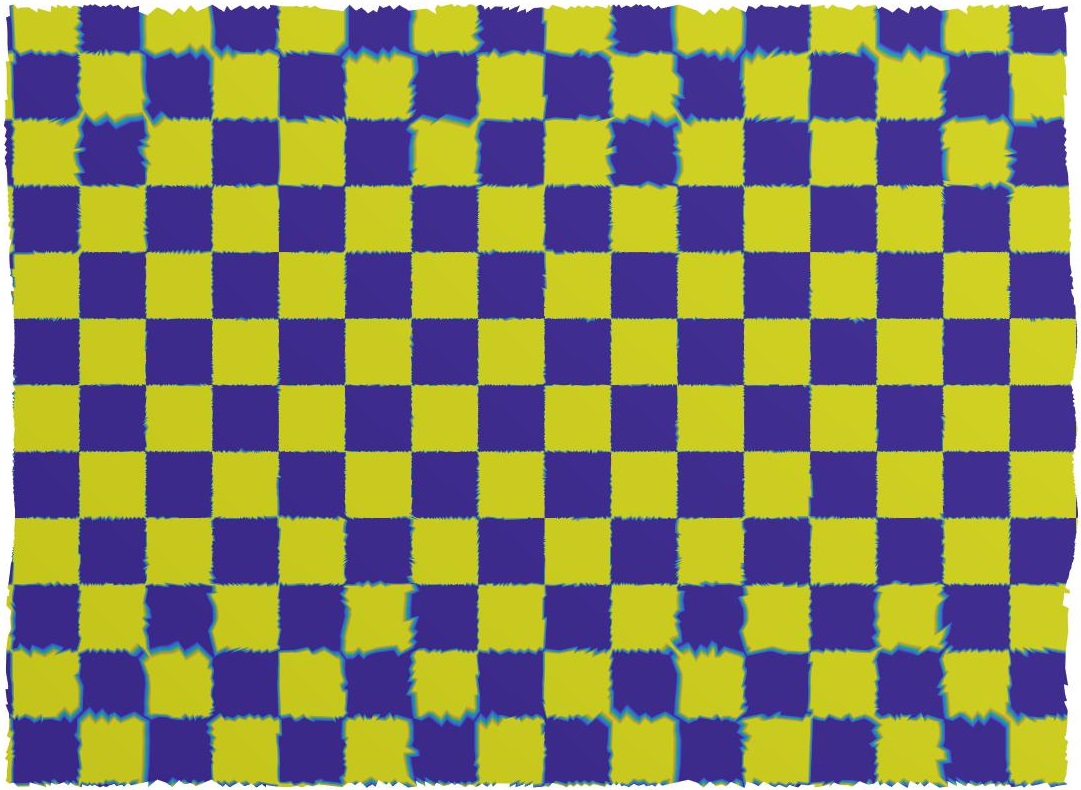}
\caption{Conformal maps to flat tori}
\end{figure}
\section{Computing Harmonic and Conformal Maps to Hyperbolic Surfaces}
This section solves the following problem on computing harmonic and conformal maps.
\begin{problem}
Assume that $M$ is a genus 1 closed smooth surface embedded in $\mathbb R^3$, and we are given
a point cloud approximation $P$ of $M$.
\begin{enumerate}[label=(\arabic*)]
\item How to compute a harmonic map
$$
f:M\rightarrow\mathbb H^2/\Gamma
$$
for a given closed hyperbolic surface $\mathbb H^2/\Gamma$ homeomorphic to $M$?

\item How to find the discrete subgroup $\Gamma$ of $Isom^+(\mathbb H^2)$ and the map
$$
f:M\rightarrow\mathbb H^2/\Gamma
$$
such that $f$ is conformal?
\end{enumerate}
\end{problem}
We will first discuss the harmonic maps in Section 7.1, and then discuss the conformal maps in Section 7.2.
\subsection{Algorithm for Harmonic Maps}

First we construct a lattice approximation $G=(V,E)$ as in Section 3, and then need to

(1) compute the discrete harmonic map $f$ from $G$ to $\mathbb H^2/\Gamma$, or equivalently, to a fundamental domain in the unit disk representation $\mathbb H^2=\{|z|<1\}$, and

(2) extend $f$ to $P$ by trilinear interpolation, and then $f|_P$ represents
an approximation of the desired harmonic map.

Step (2) is kind of clear and simple, so let us discuss the details of step (1).
First we need to "cut" the lattice $G$ along $2\cdot genus(M)$ "loops" $\gamma_i$, so that we can represent a discrete map from $G$ to a hyperbolic surface as a map from $V$ to a planar domain in the unit disk. Secondly we fix the homotopy class of discrete maps $f:G\rightarrow\mathbb H^2/\Gamma$, by requiring that each loop $\gamma_i$ corresponds to a certain deck transformation $\alpha_i$ in $\Gamma$.

Now computing the discrete harmonic map in the fixed homotopy class amounts to solve a system of nonlinear equations. For a vertex $i$ that is away from any loop $\gamma_i$, the balanced condition (1) gives us that
$$
\sum_{j:j\sim i}d_{\mathbb H^2}(f(i),f(j))\cdot\vec e(f(i),f(j))=0
$$
where $\vec e(x,y)$ denotes the unit tangent vector in $T_x\mathbb H^2$ indicating the direction to point $y$.
Things become subtle when the vertex $i$ is near some loop $\gamma_i$. Assume the edge $ij$ pass through the loop $\gamma_k$, and $l_{ij}(f)$ is small. Then
$$
f(i)\approx\alpha_{ij}(f(j)),
$$
where $\alpha_{ij}$ is a deck transformation in $\Gamma$, which is determined by

\begin{enumerate}[label=(\arabic*)]
\item the choice of the loops $\gamma_i$, and

\item the choice of $\alpha_i$'s, and

\item the relative position between $i$ and $\gamma_k$.
\end{enumerate}
Properly choosing $\alpha_i$'s and computing $\alpha_{ij}$'s are involved, particularly for higher genus surfaces. One need to use polygonal representations for higher genus surfaces, and parameterizations of the deformation spaces of $\Gamma$. See \cite{maskit1999new}\cite{maskit2001matrices} for references. Once we computed all the $\alpha_{ij}$'s, it remains to solve the following system of nonlinear equations.
\begin{equation}
\sum_{j:j\sim i}d_{\mathbb H^2}(f(i),\alpha_{ij}(f(j)))\cdot\vec e(f(i),\alpha_{ij}(f(j)))=0\quad \forall i.
\end{equation}
The existence and uniqueness of the solution is guaranteed by Theorem 2.10.
Notice that equation (10) indicates that $f(i)$ is the hyperbolic center of mass of $\alpha_{ij}(f(j))$'s where $j$ goes over all the neighbors of $i$. So one can use the so-called center of mass method to solve equation (10). This is an iterating method where in each iteration $f(i)$ is replaced by the mass-center of its neighbors $\alpha_{ij}(f(j))$, i.e., the $(n+1)$-th function $f_{n+1}(i)$ is determined by
\begin{equation*}
\sum_{j:j\sim i}d_{\mathbb H^2}(f_{n+1}(i),\alpha_{ij}(f_n(j)))\cdot\vec e(f_{n+1}(i),\alpha_{ij}(f_n(j)))=0\quad \forall i.
\end{equation*}
It has been proved by Gaster-Loustau-Monsaingeon
\cite{gaster2018computing} that the center of mass method will converge to the unique discrete harmonic map.
But
the disadvantage is that there is no direct way to compute a hyperbolic mass-center and this iteration is slow. In our numerical experiments, we use the so-called cosh-center of mass method introduced also in Gaster-Loustau-Monsaingeon's paper
\cite{gaster2018computing}. Instead of minimizing the exact discrete Dirichlet energy, we minimize
$$
{\mathcal E}_0[f]=\sum_{ij\in E}(\cosh l_{ij}(f)-1)\approx\sum_{ij\in E}\frac{l_{ij}(f)^2}{2}
$$
which is a good approximation of the discrete Dirichlet energy if all the edges lengths $l_{ij}$ are small. The balanced condition for critical points of the new energy ${\mathcal E}_0$ is
\begin{equation*}
\sum_{j:j\sim i}\sinh[d_{\mathbb H^2}(f(i),\alpha_{ij}(f(j)))]\cdot\vec e(f(i),\alpha_{ij}(f(j)))=0\quad \forall i,
\end{equation*}
and the induced iterating formula is
\begin{equation}
\sum_{j:j\sim i}\sinh[d_{\mathbb H^2}(f_{n+1}(i),\alpha_{ij}(f_n(j)))]\cdot\vec e(f_{n+1}(i),\alpha_{ij}(f_n(j)))=0\quad \forall i.
\end{equation}

The miracle is that such $\cosh$-center of mass defined as in equation (11) can be computed directly. Using the hyperboloid model
$$
\mathbb H^2_h=\{(x,y,z)\in\mathbb R^3:z^2-x^2-y^2=1\},
$$
the $\cosh$-center of points $p_1,...,p_n\in\mathbb H^2_h$ is shown, in \cite{gaster2018computing}, to be
$$
\pi_H\left(\frac{p_1+...+p_n}{n}\right)
$$
where $(p_1+...+p_n)/{n}$ is computed as vectors in $\mathbb R^3$, and $\pi_H$ is the radial projection to $\mathbb H^2_h$.

In a summary we have the following \textbf{Algorithm 6}.

\begin{algorithm}
\caption{Harmonic Maps to Hyperbolic Surfaces}

\textbf{Input:}  A point cloud $P$ in $\mathbb R^3$ as an approximation of a surface $M$, and a hyperbolic surface $\mathbb H^2/\Gamma$.
\\
\textbf{Output:}  A harmonic map $f$ from $P$ to the hyperbolic surface $\mathbb H^2/\Gamma$.

\begin{algorithmic}[1]
\State  Choose parameters $\epsilon$ and $n$, and then compute the lattice approximation $G=(V,E)$.

\State Choose $2\cdot genus(M)$ "loops" in $G$ for the cut.

\State Assign a proper deck transformation $\alpha_i\in\Gamma$ to each loop $\gamma_i$.

\State Compute the correction transformations $\alpha_{ij}$.

\State Do the $\cosh$-center of mass iterations (11), with the trivial initial map $f_0\equiv0$.

\State Extend $f$ to $P$ by trilinear interpolation in the unit disk.

\State \textbf{Return} $f|_P$.

\end{algorithmic}
\end{algorithm}

\subsection{Algorithm for Conformal Maps}
Assume that $f_\Gamma$ is the discrete harmonic map from the lattice approximation $G=(V,E)$ to the hyperbolic surface $\mathbb H^2/\Gamma$.
Inspired by part (a) of Theorem 2.5, we compute the conformal map by minimizing
$$
\frac{\mathcal E(f_\Gamma)}{Area(\mathbb H^2/\Gamma)}=\frac{\mathcal E(f_\Gamma)}{-2\pi\chi(\mathbb H^2/\Gamma)}=
\frac{\mathcal E(f_\Gamma)}{-2\pi\chi(M)}
$$
over the deformation space of $\Gamma$. If $\bar\Gamma$ is the minimizer, then
$
f_{\bar\Gamma}:V \rightarrow\mathbb H^2/{\bar\Gamma}
$
would approximate a conformal map.
Our numerical experiments used Maskit's nice parameterization of the deformation space of $\Gamma$ for genus 2 surfaces \cite{maskit1999new}. Maskit also proposed methods of parameterizing for any genus hyperbolic surfaces using Fenchel-Nielsen coordinates \cite{maskit2001matrices}.

Our method is summarized as \textbf{Algorithm 7}.
\begin{algorithm}
\caption{Conformal Maps to Hyperbolic Surfaces}

\textbf{Input:}  A point cloud $P$ in $\mathbb R^3$ as an approximation of a surface $M$.

\textbf{Output:} A hyperbolic surface $\mathbb H^2/\Gamma$ and a conformal map $f$ from $P$ to the hyperbolic surface $\mathbb H^2/\Gamma$.

\begin{algorithmic}[1]
\State  Choose parameters $\epsilon$ and $n$, and then compute the lattice approximation $G=(V,E)$.

\State Choose $2\cdot genus(M)$ "loops" in $G$ for the cut.

\State For a given hyperbolic surface $\mathbb H^2/\Gamma$ homeomorphic to $M$, assign a proper deck transformation $\alpha_i\in\Gamma$ to each loop $\gamma_i$.

\State Compute the correction transformations $\alpha_{ij}$.

\State Do the $\cosh$-center of mass iterations (11), with the trivial initial map $f_0\equiv0$.

\State Compute the discrete Dirichlet energy $\mathcal E(f_\Gamma)$.

\State Compute the minimizer $\bar \Gamma$ of $\mathcal E(f_\Gamma)$, and $f_{\bar \Gamma}$.

\State Extend $f_{\bar\Gamma}$ to $P$ by trilinear interpolation in the unit disk.

\State \textbf{Return} $f_{\bar\Gamma}|_P$.

\end{algorithmic}
\end{algorithm}
\subsection{Numerical Examples}
Here are numerical examples for harmonic maps and conformal maps on genus 2 surfaces.
\begin{figure}[H]
\centering
\includegraphics[width=0.25\textwidth]{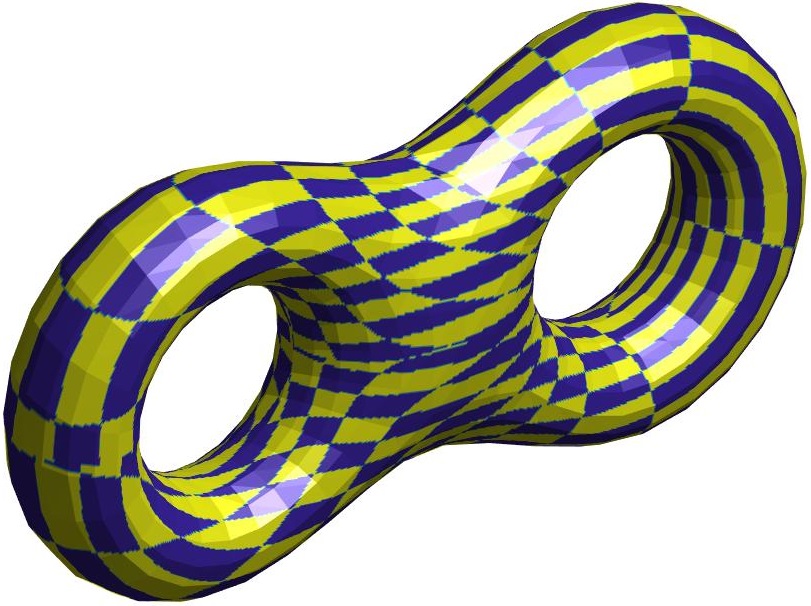}
\includegraphics[width=0.25\textwidth]{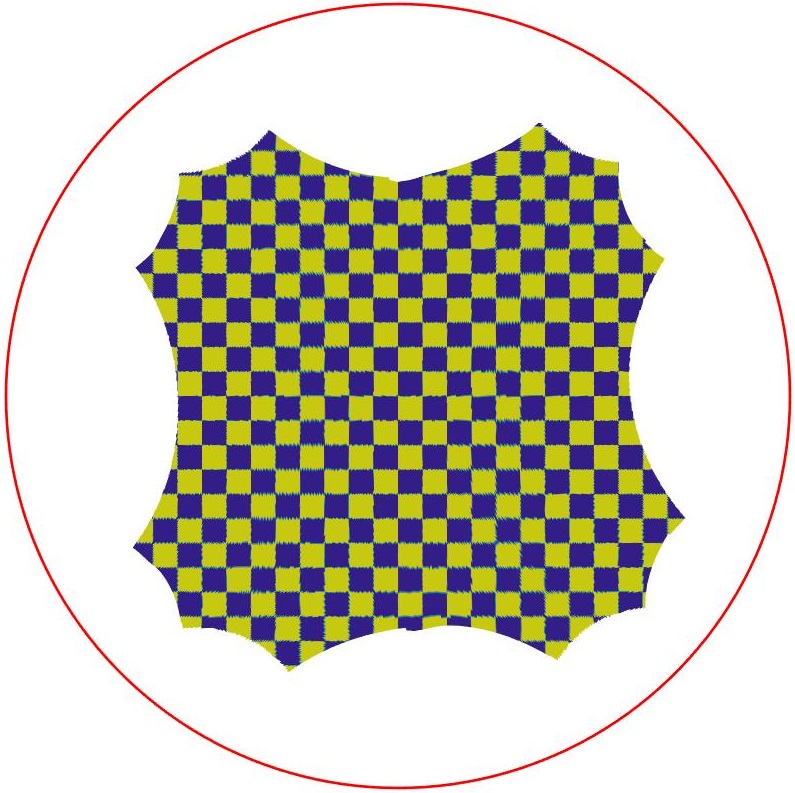}
\quad
\includegraphics[width=0.2\textwidth]{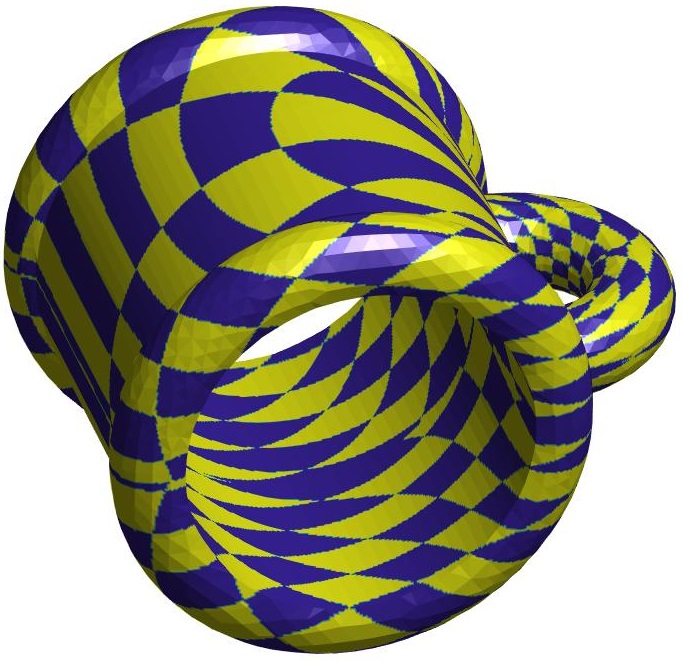}
\includegraphics[width=0.25\textwidth]{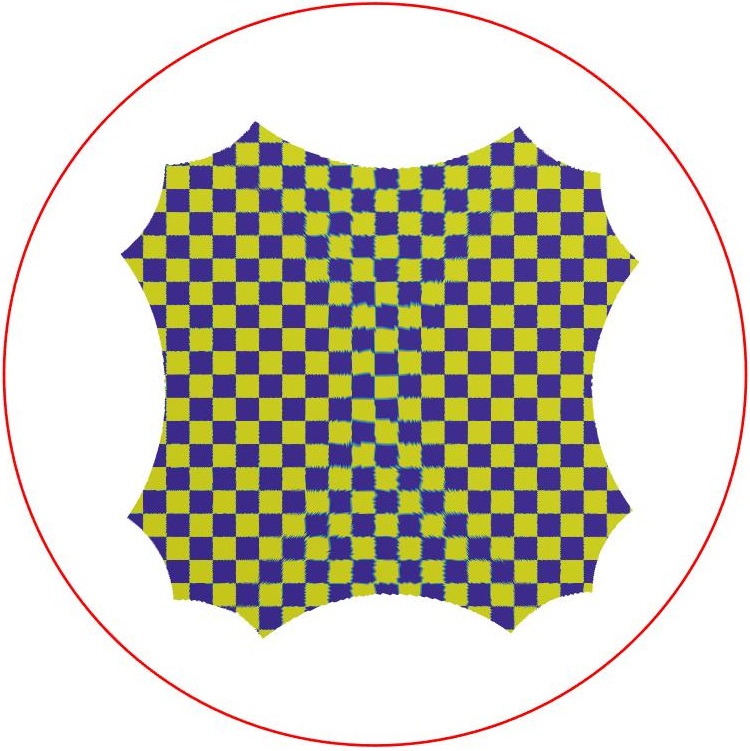}
\caption{Harmonic maps to hyperbolic surfaces}
\end{figure}

\begin{figure}[H]
\centering
\includegraphics[width=0.25\textwidth]{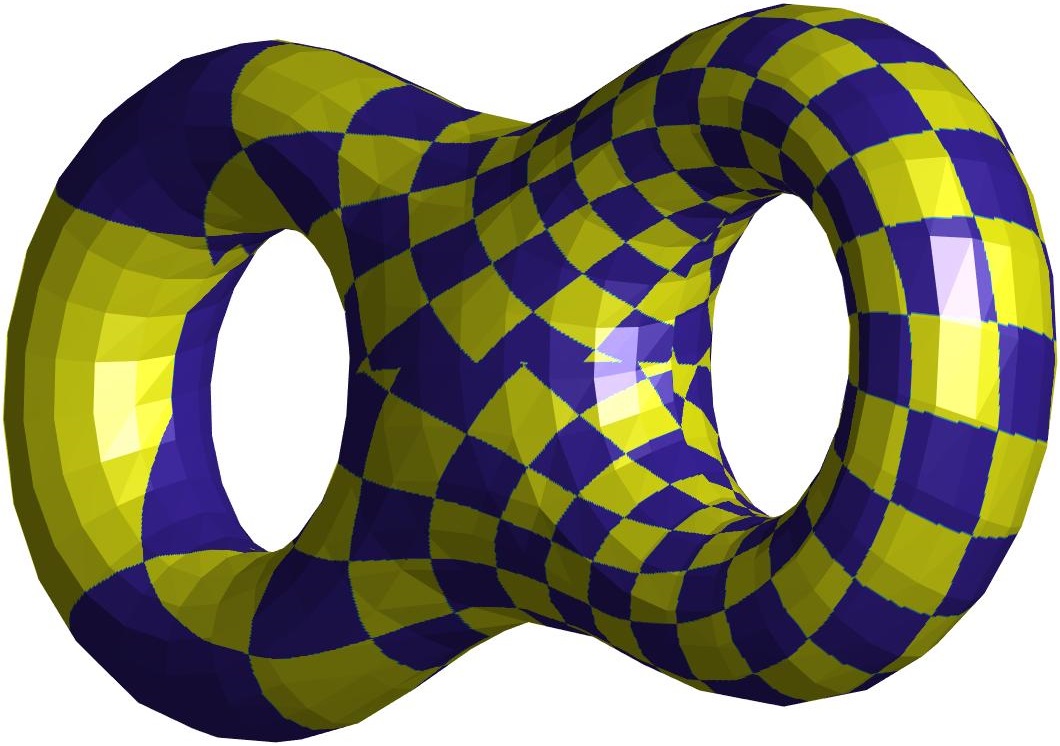}
\includegraphics[width=0.25\textwidth]{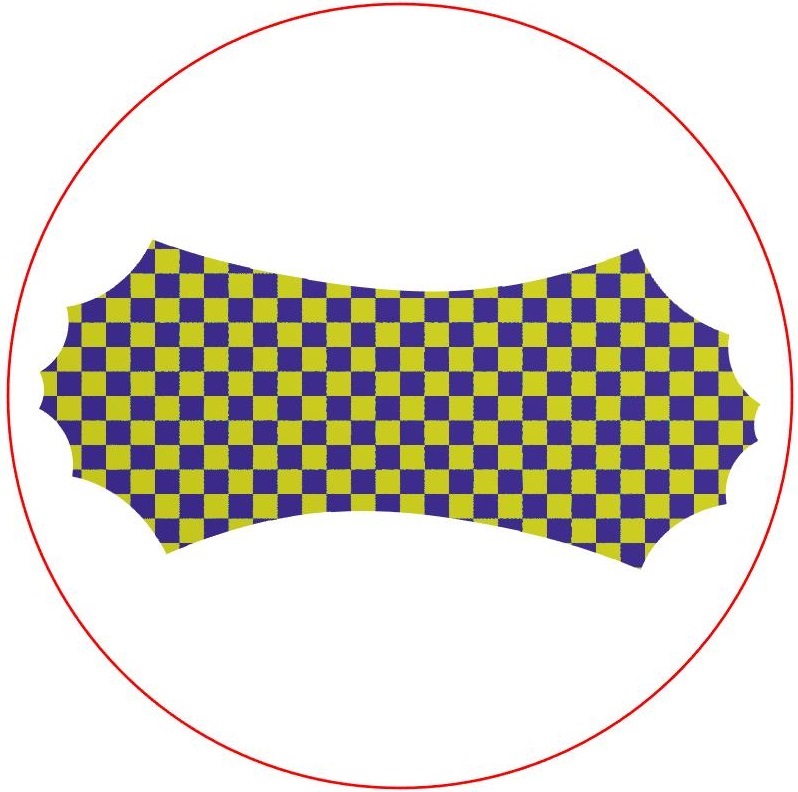}
\quad
\includegraphics[width=0.2\textwidth]{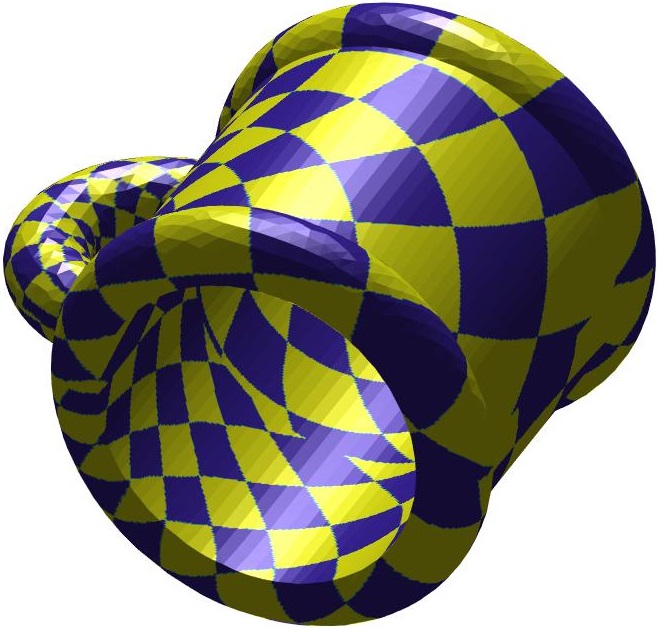}
\includegraphics[width=0.25\textwidth]{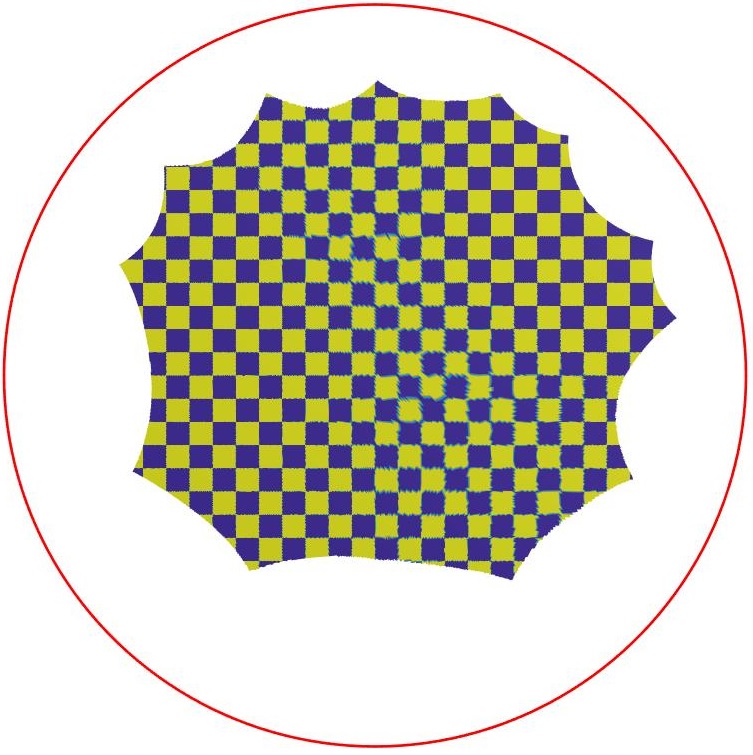}
\caption{Conformal maps to hyperbolic surfaces}
\end{figure}

\bibliographystyle{unsrt}
\bibliography{kc2gkc2}
\end{document}